\def\vers{Aug.~26, 2009, v.3}
\magnification=1200
\hsize=6.5truein
\vsize=8.9truein
\font\bigfont=cmr10 at 14pt
\font\mfont=cmr9
\font\sfont=cmr8
\font\mbfont=cmbx9
\font\sifont=cmti8
\def\Mustata{Musta\c{t}\v{a}}
\def\1{\hskip1pt}
\def\bs{\bigskip}
\def\ms{\medskip}
\def\ss{\smallskip}
\def\nin{\noindent}
\def\msum{\h{$\sum$}}
\def\mcap{\h{$\bigcap$}}
\def\mcup{\h{$\bigcup$}}
\def\mprod{\h{$\prod$}}
\def\mopl{\h{$\bigoplus$}}
\def\mwedge{\h{$\bigwedge$}}
\def\scirc{\,\raise.2ex\h{${\scriptstyle\circ}$}\,}
\def\ssb{\raise.2ex\h{${\scriptscriptstyle\bullet}$}}
\def\ssn{\,\rlap{\raise1pt\hbox{$\subset$}}{\raise-4pt
\hbox{$\1\scriptstyle\neq$}}\,}
\def\don{\rlap{$\downarrow$}\raise-2pt\hbox{$\downarrow$}}
\def\a{\alpha}
\def\B{{\cal B}}
\def\C{{\bf C}}
\def\d{\partial}
\def\D{{\cal D}}
\def\e{\widetilde{e}}
\def\ee{\widetilde{e}^{\,\prime}}
\def\E{{\cal E}}
\def\F{{\cal F}}
\def\f{\widetilde{f}}
\def\g{\gamma}
\def\G{{\cal G}}
\def\h{\hbox}
\def\H{\widetilde{H}}
\def\I{{\cal I}}
\def\J{{\cal J}}
\def\K{{\cal K}}
\def\l{\lambda}
\def\L{\Lambda}
\def\N{{\bf N}}
\def\O{{\cal O}}
\def\P{{\bf P}}
\def\q{\quad}
\def\Q{{\bf Q}}
\def\R{{\bf R}}
\def\S{{\cal S}}
\def\SS{\bar{\cal S}}
\def\X{\widetilde{X}}
\def\Y{\widetilde{Y}}
\def\Z{\widetilde{Z}}
\def\ZZ{{\bf Z}}
\def\tcS{\widetilde{\cal S}}
\def\Gr{\h{\rm Gr}}
\def\DR{\h{\rm DR}}
\def\JC{\h{\rm JC}}
\def\Ker{\h{\rm Ker}}
\def\Coker{\h{\rm Coker}}
\def\Sp{\h{\rm Sp}}
\def\Spec{{\cal S}pec}
\def\rank{{\rm rank}}
\def\nnc{{\rm nnc}}
\def\nrnc{{\rm nrnc}}
\def\red{{\rm red}}
\def\codim{\h{\rm codim}}
\def\mult{\h{\rm mult}}
\def\simto{\buildrel\sim\over\longrightarrow}
\def\lc{\raise1pt\hbox{$\lceil$}}
\def\rc{\raise1pt\hbox{$\rceil$}}
\def\lf{\lfloor}
\def\rf{\rfloor}
\def\lb{\bigl[}
\def\rb{\bigr]}
\def\lkr{\langle k\rangle}
\h{}
\vskip 1cm
\centerline{\bigfont Jumping coefficients and spectrum}

\ss
\centerline{\bigfont of a hyperplane arrangement}

\bs
\centerline{Nero Budur and Morihiko Saito}

\bs\ms
{\narrower\nin
{\mbfont Abstract.} {\mfont
In an earlier version of this paper written by the second named
author, we showed that the jumping coefficients of a hyperplane
arrangement depend only on the combinatorial data of the
arrangement as conjectured by \Mustata.
For this we proved a similar assertion on the spectrum.
After this first proof was written, the first named author found a
more conceptual proof using the Hirzebruch-Riemann-Roch theorem
where the assertion on the jumping numbers was proved without
reducing to that for the spectrum.
In this paper we improve these methods and show that the jumping
numbers and the spectrum are calculable in low dimensions without
using a computer.
In the reduced case we show that these depend only on fewer
combinatorial data, and give completely explicit combinatorial
formulas for the jumping coefficients and (part of) the spectrum
in the case the ambient dimension is 3 or 4.
We also give an analogue of \Mustata's formula for the spectrum.}
\par}

\bs\bs
\centerline{\bf Introduction}
\footnote{}{{\sifont Date\1}\sfont:\ \vers}
\footnote{}{\sfont The first author is supported by the NSF grant
DMS-0700360.}
\footnote{}{\sfont The second author is partially supported by
Kakenhi 19540023.}

\bs
This paper combines and improves two unpublished preprints: [29]
which gave the first proof of Theorem~1 below, and [5] which gave a
formula for the spectrum using the Hirzebruch-Riemann-Roch theorem
[19] together with the combinatorial description of the cohomology
ring of the wonderful model [8].

Let $D$ be a hyperplane arrangement in $X=\C^n$ with $D_i\,(i\in\L)$
the irreducible components.
In this paper we assume $D$ is central and essential
(i.e.\ $0\in D_i$ for any $i$, and
$\bigcap_iD_i=\{0\}$) since the multiplier ideals and the spectrum
are defined locally.
We also assume that $D$ is not a divisor with normal crossings
(i.e.\ $\deg D_{\red}>n$).
However, we do {\it not\1} assume $D$ is {\it reduced}.

For a positive rational number $\a$, the multiplier ideal
$\J(X,\a D)$ is a coherent ideal of the structure sheaf
$\O_{X}$ defined by the local integrability of $|g|^{2}/
|f|^{2\a}$ for $g\in\O_{X}$, where $f$ is a defining polynomial of
$D$, see [21], [23].
This gives a decreasing sequence of ideals $\J(X,\a D)$ which is
locally constant outside a locally finite subset $\JC(D)$ of $\Q$.
The elements of $\JC(D)$ are called the {\it jumping coefficients}.
Here we may restrict to the intersection with the interval $(0,1)$
since $1\in\JC(D)$ and for $\a>0$, we have $\a\in\JC(D)$ if and
only if $\a+1\in\JC(D)$.
A formula for the multiplier ideals of a hyperplane arrangement
$D$ was obtained by {\Mustata} [22]. It was conjectured there that
$\JC(D)$ depends only on the combinatorial data of a hyperplane
arrangement.

For simplicity, consider the case where $n=3$, $D$ is reduced, and
$\mult_xD\le 3$ for any $x\ne 0$.
Let $D^{\nnc}\subset D$ denote the complement of the subset
consisting of normal crossing singularities.
This is the union of the lines of multiplicity 3, and corresponds
to a finite subset $\P(D^{\nnc})\subset\P^2$.
Let $\C[x,y,z]_i$ denote the vector space of homogeneous polynomials
of degree $i$. Set $d=\deg D$.
In this case {\Mustata}'s formula implies the following
(see [22], Cor.~2.1):
If $\a=j/d\in({2\over 3},1)$ with $j\in\ZZ$, then
$$\a\in\JC(D)\iff\exists\,g\in\C[x,y,z]_{j-3}\setminus\{0\}
\,\,\,\h{with}\,\,\,g^{-1}(0)\supset\P(D^{\nnc}).\leqno(0.1)$$
The last condition may a priori depend on the position of
$\P(D^{\nnc})$, and the conjecture is rather nontrivial.
It turns out, however, that the points of $\P(D^{\nnc})$ are always
in a generic position as far as the above condition is concerned,
see a remark on the surjectivity of (0.2) after Theorem~5 below.
Note that there is an example such that ${j\over d}\notin
\JC(D)$, see Example~(3.6) below.
(This implies a negative answer to [22], Question~2.7
even in the case $j=5$.)

Using an inductive argument, we gave in [29] the first proof of
the following.

\ms\nin
{\bf Theorem~1.} {\it The jumping coefficients and the spectrum of
a hyperplane arrangement depend only on the combinatorial data.}

\ms
The spectrum $\Sp(f)$ of a hypersurface singularity was defined by
J.~Steenbrink [31] using the monodromy and the Hodge filtration on
the Milnor cohomology where $f$ is a function defining locally $D$.
The spectrum $\Sp(f)$ is a fractional polynomial
$\sum_{\a}n_{f,\a}t^{\a}$ where $n_{f,\a}=0$ for $\a\notin\Q\cap
(0,n)$ (see [6], Prop.~5.2), and $\a$ is called a nontrivial
exponent if $n_{f,\a}\ne 0$.
It was shown in [3] that, for $\a\in(0,1)$, it is a jumping
coefficient if it is a nontrivial exponent,
and the converse holds in the isolated singularity case.
Using [6], we can show for any holomorphic function $f$ that the
converse also holds if
$\a$ is an isolated jumping coefficient (i.e.\ if $\a$ is not a
jumping coefficient for $D\setminus\{0\}$),
see Proposition~(4.2) below.

Returning to the case of an essential central hyperplane arrangement,
let $g_i$ be the defining linear function of $D_i$ with $m_i$ the
multiplicity.
Let $V$ be an intersection of $D_i$ which is called an {\it edge}.
Set $f_{X/V}=\prod_{D_i\supset V}g_i^{m_i}$.
This is viewed as a function on $X/V$, and defines a hyperplane
arrangement $D_V\subset X/V$.
Let $D^{\nrnc}\subset D$ denote the complement of the subset
consisting of reduced normal crossing singularities.
By Proposition~(4.2) we then get the following.

\ms\nin
{\bf Proposition~1.} {\it For $\a\in(0,1)$, we have $\a\in\JC(D)$
if and only if there is an edge $V\subset D^{\nrnc}$ with
$n_{f_{X/V},\a}\ne 0$.}

\ms
So the assertion of Theorem~1 for the jumping coefficients is
reduced to the corresponding assertion for the spectrum since the
combinatorial data of $D_V\subset X/V$ are obtained by restricting
those of $D$.
The spectrum of a hyperplane arrangement is calculated by using the
canonical embedded resolution together with the filtered logarithmic
complexes associated to certain local systems as in [16].
The assertion is then reduced to the calculation of the restriction
of the de Rham complex to the exceptional divisors of
the canonical embedded resolution where we need some arguments
as in [8].
Note that the moduli space of hyperplane arrangements with fixed
combinatorial data is not necessarily connected as shown in [25]
as a corollary of a deep theorem (see also [27], 5.7 for a direct
argument), and hence we cannot prove Theorem~1 by simply using a
deformation argument.

After an earlier version of the above proof of Theorem~1 due to
the second named author [29] was written, the first named author
([4], [5]) found a more conceptual proof using the
Hirzebruch-Riemann-Roch theorem and a combinatorial description [8]
of the cohomology ring of the canonical embedded resolution of the
corresponding hyperplane arrangement in $\P^{n-1}$.
The new proof implies formulas for the jumping coefficients and
the spectrum in terms of the combinatorial data where the
assertion for the jumping numbers is proved without using
Proposition~1 (although it needs a compactification of $X=\C^n$
so that the calculation becomes more complicated than the proof
using Proposition~1).
It is also possible to write down a formula for the spectrum as in
[5] by summarizing the arguments in (5.3--4) in this paper.
(Note that [4] did not deal with the whole spectrum.)
Stimulated by this new proof, there was an improvement of the
inductive argument.
Using these we can prove Theorem~2 below.

Let $\S(D)^{\nnc}$ denote the set of edges $V$ contained in
$D^{\nnc}$.
Let $\subset$, $\mu(V)$, and $\g(V)$ respectively denote
the inclusion relation, multiplicity of $D$ along $V$, and the
codimension of $V$, see (1.1) below.
Then we have the following.

\ms\nin
{\bf Theorem~2.} {\it Assume $D$ is reduced. Then the jumping
coefficients and the coefficients $n_{f,\a}$ of the spectrum for
$\a\in(0,1]\cup(n-1,n)$ depend only on the weak combinatorial
equivalence class, i.e.\ on the set $\S(D)^{\nnc}$ together with
the combinatorial data $\subset$, $\mu$, $\g$.
For $\a\in(1,n-1]$, $n_{f,\a}$ depends only on $\S(D)^{\nnc}$
together with $\subset$, $\mu$, $\g$ and also the subsets
$\S^{D_i}:=\{V\in\S(D)^{\nnc}\mid V\subset D_i\}\,(i\in\L)$.}

\ms
The weak combinatorial equivalence is strictly weaker than the
usual combinatorial equivalence.
Indeed, if $n=3$, the former is determined
only by $d=\deg D$ and
$$\nu^{(2)}_m=\#\{y\in\P(D)\mid\mult_y\P(D)=m\}\q\h{for}\,\,\,
m\ge 3,$$
where $\P(D)\subset\P^2$, see also (1.1.3).
For instance, in the case where $d=7$, $\nu^{(2)}_3=3$, and
$\nu^{(2)}_m=0\,(m>3)$,
there are two possibilities of combinatorial data depending on
whether the three points of multiplicity 3 are on a same line or
not.

In this paper we also show that the jumping coefficients and
the spectrum of hyperplane arrangements are calculable in low
dimensions without using a computer as in [4], [5]
(although the formula is rather complicated).
This was partly made possible by restricting the centers of the
blow-ups to those contained in $D^{\nnc}$.
For instance we have the following.

\ms\nin
{\bf Theorem~3.} {\it Assume $D$ reduced and $n=3$.
Then $n_{f,\a}=0$ if $\a d\notin\ZZ$, and we have for
$\a={i\over d}\in(0,1]$ with $i\in[1,d]\cap\ZZ$
$$\eqalign{n_{f,\a}&=\h{${i-1\choose 2}-\sum_m
\nu^{(2)}_m{\lceil im/d\rceil-1\choose 2}$},\cr
n_{f,\a+1}&=(i-1)(d-i-1)-\msum_m\,\nu^{(2)}_m
\bigl(\lc im/d\rc-1\bigr)\bigl(m-\lc im/d\rc\bigr),\cr
n_{f,\a+2}&=\h{${d-i-1\choose 2}-\sum_m\nu^{(2)}_m
{m-\lceil im/d\rceil\choose 2}-\delta_{i,d}$},}$$
where $\lc\beta\rc:=\min\{k\in\ZZ\mid k\ge\beta\}$,
and $\delta_{i,d}=1$ if $i=d$ and $0$ otherwise.}

\ms\nin
{\bf Theorem~4.} {\it Assume $D$ reduced and $n=4$.
Let $\nu^{(2)}_m$, $\nu^{(3)}_{m'}$, $\nu^{(2,3)}_{m,m'}$ be as
in $(1.1.3)$ below.
Then $n_{f,\a}=0$ for $\a d\notin\ZZ$, and we have for
$\a={i\over d}\in(0,1]$ with $i\in[1,d]\cap\ZZ$}
$$\eqalign{n_{f,\a}&=\h{${i-1\choose 3}-\msum_{m,m'}\,
\nu^{(2,3)}_{m,m'}\Bigl(2{\lceil im/d\rceil-1\choose 3}-
{\lceil im/d\rceil-1\choose 2}\big(\lc im'/d\rc-3\bigr)\Bigr)$}\cr
&+\h{$\msum_m\,\,\nu^{(2)}_m\Bigl(2{\lceil im/d\rceil-1
\choose 3}-(i-3){\lceil im/d\rceil-1\choose 2}\Bigr)
-\msum_{m'}\,\nu^{(3)}_{m'}{\lceil im'/d\rceil-1\choose 3}$}.}$$

\ms
The formula is similar for $n=4$ and $\a\in(3,4)$.
However, it requires some more combinatorial data, and is more
complicated for $n=4$ and $\a\in(1,3]$.
Those are left to the reader as exercises.
Note that the formula for $n=3$ and $\a\in(0,1]$ is quite similar
to a formula for the Hodge number in [15], Th.~6.
As for the jumping coefficients, it is well-known that
$\JC(D)\cap(0,1)=\{i/d\mid i\in\ZZ\cap[2,d)\}$ with $d=\deg D$ if
$n=2$ and $D$ is reduced.
Combined with Proposition~1, Theorems~3--4 then imply the
following.

\ms\nin
{\bf Corollary~1.} {\it Assume $D$ reduced and $n=3$.
Then $\a\in(0,1)$ is a jumping coefficient of $D$ if and only if
there is $m\in\N\cap[3,\infty)$ such that $m\a\in\ZZ\cap[2,m)$ and
$\nu^{(2)}_m\ne 0$ or there is $i\in\ZZ\cap[3,d)$ such that
$\a={i\over d}$ and $n_{f,\a}\ne 0$ in Theorem~$3$.}

\ms\nin
{\bf Corollary~2.} {\it Assume $D$ reduced and $n=4$.
In the notation of $(1.1.3)$, $\a\in(0,1)$ is a jumping coefficient
of $D$ if and only if there is $m\in\N\cap[3,\infty)$ such that
$m\a\in\ZZ\cap[2,m)$ and $\nu^{(2)}_m\ne 0$, or
there is $V\in\S(D)^{\nnc}$ together with $i\in\ZZ\cap[3,\mu(V))$
such that $\codim\,V=3$, $\a=i/\mu(V)$ and $n_{f_{X/V},\a}\ne 0$
in Theorem~$3$, or
there is $i\in\ZZ\cap[4,d)$ such that $\a=i/d$ and $n_{f,\a}\ne 0$
in Theorem~$4$.
Here $\mu(V)=\mult_VD$.}

\ms
We have a formula for the spectrum analogous to {\Mustata}'s
formula [22] as follows.

\ms\nin
{\bf Theorem~5.} {\it With the notation of $(1.1)$ below,
let $\I_V\subset\C[X]$ denote the reduced ideal of $V\in
\S':=\S(D)^{\nrnc}$.
For $\a={j\over d}\in(0,1]$ with $j\in[1,d]\cap\ZZ$ we have}
$$n_{f,\a}=\dim\Bigl(\mcap_{V\in\S'\setminus\{0\}}\,
\I_V^{\lceil\a\mu(V)\rceil-\g(V)}\cap\C[X]_{j-n}\Bigr).$$

By Proposition~(4.2) below this is compatible with {\Mustata}'s
formula [22] for $\a\in(0,1)$ in the case $\mu(V)\a\notin\ZZ$ for
any $V\in\S(D)\setminus\{0\}$, see also (3.3) below.
Here $\C[X]_{j-n}$ denotes the space of homogeneous polynomials of
degree $j-n$, and the intersection with $\I_V^{e_V}$ gives a
restriction condition for $g\in\C[X]_{j-n}$ as in the right-hand
side of (0.1) where $e_V=\lceil\a\mu(V)\rceil-\g(V)$.
In the case $n=3$ and $D$ is reduced, this restriction condition is
given by the condition that $g\in{\bf m}_V^{e_V}$, where $g$ is
viewed as a section of $\O_{\P^2}(j-3)$ and ${\bf m}_V\subset
\O_{\P^2,y}$ is the maximal ideal with $y$ the closed point of
$\P^2$ corresponding to $V$.
(Here $\O_{\P^2}(j-3)$ is locally trivialized.)
Since $\dim\O_{\P^2,y}/{\bf m}_V^{e_V}={e_V+1\choose 2}$, the first
equality of Theorem~3 for $\a\in(0,1]$ implies the surjectivity of
$$H^0(\P^2,\O_{\P^2}(j-3))\to\mopl_{V\in\S'\setminus\{0\}}\,
\O_{\P^2,y}/{\bf m}_V^{e_V},\leqno(0.2)$$
which means that the above restriction conditions for
$V\in\S'\setminus\{0\}$ are always independent so that they give
the maximal restriction condition in total.
This is closely related to the non-degeneracy of the matrix after
(3.6.1).
Note also that the above surjectivity is equivalent to the
vanishing
$$H^1(\P^2,\O_{\P^2}(i-3)\otimes_{\O}\I(\a))=0,$$
where $\I(\a):=\Ker(\O_{\P^2}\to\mopl_{V\in\S'\setminus\{0\}}\,
\O_{\P^2,y}/{\bf m}_V^{e_V})$.

It is known that the jumping coefficients are closely related to
the roots of the $b$-function, see [14] and (3.5) below.
For the moment it is unclear whether an analogue of Theorem~1
holds for the $b$-function.
As for the relation with topology, note that the spectrum does
not determine each Betti number of the Milnor fiber of an
hyperplane arrangement, see [7], Ex.~5.4--5 and [12], p.~213,
Ex.~4.16.

We would like to thank M.~{\Mustata} and A.~Dimca for useful
comments concerning this paper.
We also thank the referee for useful remarks.

In Section~1, we review some facts related to hyperplane
arrangements, spectrum and jumping coefficients.
In Section~2, we essentially reproduce Section~2 of [29] on the
canonical embedded resolution of a projective hyperplane
arrangement, see also [8].
In Section~3, we prove Theorems~3 and 5.
In Section~4, we give an improved version of the first proof of
Theorem~1 together with proofs of Theorems~2--4 by induction.
In Section~5, we improve some arguments in [4], [5] (using
$D^{\nnc}$), and give proofs of Theorems~1--4 using the
Hirzebruch-Riemann-Roch theorem and the combinatorial description
of the cohomology ring of the embedded resolution.

\bs\bs
\centerline{\bf 1. Preliminaries}

\bs\nin
{\bf 1.1.~Hyperplane arrangements.}
Let $D$ be a central hyperplane arrangement in $X=\C^n$ with $D_i\,
(i\in\L)$ the irreducible components of $D$, where central
means that all the $D_i$ pass through the origin, see [24].
We also assume $D$ is essential, and is not
a divisor with normal crossings, i.e.\ $\bigcap_iD_i=\{0\}$ and
$\deg D_{\red}>n$.
We define a set of vector subspaces of $X$ by
$$\S(D)=\{\mcap_{i\in I}\,D_i\}_{I\subset\L,I\ne\emptyset},$$
where $I$ runs over the nonempty subsets of $\L$.
(Note that we may have $\mcap_{i\in I}\,D_i=\mcap_{i\in I'}\,D_i$
with $I\ne I'$.)
For $V\in\S(D)$, define
$$I(V)=\{i\in\L\mid D_i\supset V\},$$
so that $V=\mcap_{i\in I(V)}\,V_i$.
There are functions $\g,\mu,\mu_{\red}:\S(D)\to\N$ such that
for $V\in\S(D)$
$$\eqalign{\g(V)&:=\codim_XV=
\min\{|I|\,\,\big|\,\,\mcap_{i\in I}\,D_i=V\},\cr
\mu(V)&:=\h{mult}_VD=\msum_{i\in I(V)}\,\mu(D_i),\cr
\mu_{\red}(V)&:=\h{mult}_VD_{\red}=\# I(V).}
\leqno(1.1.1).$$
There is a natural order on $\S(D)$ defined by the inclusion
relation $\subset$.

Let $D^{\nnc}\subset D$ denote the complement of the subset
consisting of normal crossing singularities.
Here ``normal crossing'' means that the associated reduced variety
has normal crossings.
We define similarly $D^{\nrnc}$ by replacing ``normal crossing''
with ``reduced normal crossing'' so that $D^{\nnc}\subset D^{\nrnc}$.
(In this paper we do not assume $D^{\nnc}\ne\emptyset$.)
Set
$$\S(D)^{\nnc}=\{V\in\S(D)\mid V\subset D^{\nnc}\}\,\,\,
\h{(similarly for}\,\,\S(D)^{\nrnc}).\leqno(1.1.2)$$
In this paper $D^{\nrnc}$ is used only in Theorem~5 and
{\Mustata's} formula, see (3.3--4).

For $\S:=\S(D)^{\nnc}$, set
$$\eqalign{&\S^{D_i}=\{V\in\S\mid V\subset D_i\},\q
\S^{(i)}=\{V\in\S\mid\g(V)=i\},\cr
&\nu^{(i)}_m=\#\{V\in\S^{(i)}\mid\mu(V)=m\},\cr
&\nu^{(i,j)}_{m,m'}=\#\{(V,V')\in\S^{(i)}\times\S^{(j)}\mid
V\supset V',\,\mu(V)=m,\,\mu(V')=m'\}.}\leqno(1.1.3)$$
This is compatible with the definition of $\nu^{(2)}_m$ in
Introduction.

\ms\nin
{\bf 1.2.~Combinatorial equivalence class.}
With the above notation we call
$$\S(D),\,\subset,\,\mu$$
the (strong) combinatorial data of a hyperplane arrangement $D$.
Note that $\g$ is determined by the inclusion relation
$\subset$, and so is $\mu$ in the case $D$ is reduced.
We call
$$\S(D)^{\nnc},\,\subset,\,\mu,\,\g,$$
the {\it weak\1} combinatorial data.
We say that two hyperplane arrangements $D$ and $D'$ are
combinatorially equivalent if there is a one-to-one correspondence
between $\S(D)$ and $\S(D')$ in a compatible way with
$\subset,\mu$ (similarly for weak combinatorial equivalence).

\ms\nin
{\bf 1.3.~Milnor fiber and the covering.}
In this subsection $D\subset X:=\C^n$ is defined by a homogeneous
polynomial $f$. Set $Z=\P(D)$.
This is a subvariety of $Y:=\P^{n-1}$ defined by $f$.
There is a ramified covering
$$Y':=\Spec_{\,Y}(\mopl_{0\le k<d}\,\S^k)\buildrel\pi\over\to Y,$$
where $\S^k=\O_Y(-k)$ and $f$ induces morphisms
$\O_Y(-k-d)\to\O_Y(-k)$ defining a ring structure on
$\mopl_{0\le k<d}\,\S^k$.

Let $U=Y\setminus Z$, and $U'$ be the restriction of $Y'$ over
$U$. Then $U'$ is \'etale over $U$, and the Milnor fiber
$f^{-1}(1)$ is identified with $U'$.
Indeed, $Y'$ is identified with a section of the line bundle
over $Y$ corresponding to $\O_Y(1)$, and $U'$ is identified
with a section of its dual bundle which is isomorphic to the
blow-up of $\C^n$ at the origin. So $U'$ is identified with
the divisor on $\C^n$ defined by $f=1$.

The geometric Milnor monodromy is induced by an automorphism of
$\C^n$ defined by
$$T_g:(x_1,\dots,x_n)\mapsto(\xi x_1,\dots,\xi x_n),$$
where $\xi=\exp(2\pi\sqrt{-1}/d)$.
Note that the monodromy of the cohomology local system
associated with the Milnor fibration on a punctured disk
is given by $(T_g^*)^{-1}$.
It is well-known (see e.g.\ [7], [12]) that
$$\S^k=\Ker((T_g^*)^{-1}-\xi^k)\subset\pi_*\O_{Y'}=
\mopl_{0\le k<d}\,\S^k.\leqno(1.3.1)$$
For the convenience of the reader, we include here a proof.
Using the projective coordinates $z_0,\dots,z_n$ of
$\P^n\supset\C^n$ such that $x_i=z_i/z_0$ for $i\in[1,n]$, the
geometric monodromy is induced by
$$T_g:(z_0,z_1,\dots,z_n)\mapsto(\xi^{-1} z_0,z_1,\dots,z_n).$$
After changing the order of the coordinates $z_1,\dots,z_n$ if
necessary, let $y_i=z_i/z_n$ on $\{z_n\ne 0\}\subset\P^n$ for
$i\in[0,n-1]$.
Set $h(y_0,\dots,y_{n-1})=f(z_1,\dots,z_n)/z_n^d$.
Then $Y'\subset\P^n$ is locally defined by the equations
$$f(x_1,\dots,x_n)=1,\q f(z_1,\dots,z_n)=z_0^d,\q
h(y_0,\dots,y_{n-1})=y_0^d,$$
and the restriction of $\pi_*\O_{Y'}$ to $\C^{n-1}=\{z_n\ne 0\}
\subset\P^{n-1}$ is identified with
$$\pi_*(\O_{\C^n}/(y_0^d-h(y_0,\dots,y_{n-1}))=
\mopl_{k=0}^{d-1}\,\O_{\C^{n-1}}y_0^k,$$
where $\{z_n\ne 0\}\subset\P^n$ is identified with $\C^n$.
On the other hand, the action of $T_g^*$ on the coordinate $y_i$
is the multiplication by $\xi^{-1}$ for $i=0$, and the identity
for $i\ne 0$.
So (1.3.1) follows.

\ms\nin
{\bf 1.4.~Spectrum.}
With the notation of (1.3) the spectrum
$\Sp(f)=\msum_{\a\in\Q}\,n_{f,\a}t^{\a}$ is defined by
$$\eqalign{&n_{f,\a} =\msum_{j}\,(-1)^{j-n+1}\dim\Gr_{F}^{p}
\H^j(f^{-1}(1),\C)_{\l}\cr
&\q\q\h{with}\,\, p=\lf n-\a\rf,\,\,\l=\exp(-2\pi i\a),}$$
where $\H^j(f^{-1}(1),\C)_{\l}$ is the
$\l$-eigenspace of the reduced cohomology of $f^{-1}(1)
\subset\C^n$ for the semi-simple part of the Milnor monodromy,
and $F$ is the Hodge filtration, see [31].
Here $\lf\beta\rf:=\max\{k\in\ZZ\mid k\le\beta\}$.
By [6], Prop.~5.2, we have
$$n_{f,\a}=0\q\h{if}\q\a\notin(0,n).$$

By (1.3.1) $\S^k$ has a meromorphic connection with
regular singularities along $Z$, and hence the localization
$\S^k(*Z)$ along $Z$ is a regular holonomic $\D_Y$-module,
which is locally isomorphic to $\O_Y(*Z)h^{k/d}$ where $h$ is as
in (1.3).
(This is proved by using the equation $h=y_0^d$.)
Moreover we get by (1.3.1)
$$H^j(U,\DR(\S^k|_U))=H^j(f^{-1}(1),\C)_{\l}\q\h{with}\,\,\,
\l=\exp(2\pi ik/d).\leqno(1.4.1)$$
Since $\S^k$ has rank 1, this implies
$$\msum_{j\in\ZZ}\,(-1)^j\dim H^j(f^{-1}(1),\C)_{\l}=\chi(U).
\leqno(1.4.2)$$

Let $\rho:\Y\to Y$ be an embedded resolution of $Z$ inducing an
isomorphism over $Y\setminus Z$.
We have a divisor with normal crossings on $\Y$
$$\Z:=\rho^*Z=Z'+Z''\q\h{with}\,\,\,Z'=\msum_{j\in J'}\,m_jE_j,
\,\,\,Z''=\msum_{j\in J''}\,m_jE_j,$$
where $Z'$ is the proper transform of $Z$, $Z''$ is the exceptional
divisor, and the $E_j$ are the irreducible components with
multiplicity $m_j$. Set $J=J'\cup J''$.

Let $\tcS^k$ be the Deligne extension of $\S^k|_U$ over
$\Y$ such that the eigenvalues of the residue of the connection are
contained in $[0,1)$, see [9].
Let $\H$ be the total (or proper) transform of a general hyperplane
$H$ of $Y$.
Using the pull-back of the above local form $\O_Y(*Z)h^{k/d}$,
we get
$$\tcS^k=\O_{\Y}(-k\H+\msum_{j\in J}\,\lf km_j/d\rf E_j),$$
since $\O_Y(Z)=\O_Y(dH)$.
Indeed, the eigenvalue of the residue of the connection along $E_j$
is
$$km_j/d-\lf km_j/d\rf.$$
Note that the above summation is taken over $J''$ in case $D$ is
reduced since $\lf km_j/d\rf=0$ if $m_j=1$ and $k\in[0,d)$.
It is known (see e.g.\ [16]) that the associated filtered
logarithmic complex together with the filtration $\sigma$
(see [10]) calculates the Hodge filtration on the cohomology.
It coincides with the Hodge filtration obtained from the
theory of mixed Hodge modules.
(This is shown by using e.g.\ [26], 3.11.)
So we get
$$\Gr_F^pH^{p+q}(f^{-1}(1),\C)_{\l}
=H^q(\Y,\Omega_{\Y}^{p}(\log\Z)\otimes_{\O}\O_{\Y}
(-k\H+\msum_j\,\lf km_j/d\rf E_j)).\leqno(1.4.3)$$
Since $f^{-1}(1)$ is affine and $(n-1)$-dimensional,
we have for $q>n-1-p$
$$H^q(\Y,\Omega_{\Y}^{p}(\log\Z)\otimes_{\O}\O_{\Y}
(-k\H+\msum_{j\in J}\,\lf km_j/d\rf E_j))=0.\leqno(1.4.4)$$
This is closely related to the Kodaira-Nakano type vanishing
theorem in [16].
We also have
$$\eqalign{\Omega_{\Y}^{n-1}(\log\Z)
&=\O_{\Y}(-n\H+\msum_{j\in J}\,c_jE_j),\cr
&=\O_{\Y}((d-n)\H+\msum_{j\in J}\,(c_j-m_j)E_j),}\leqno(1.4.5)$$
where $c_j$ is the codimension of the center of the blow-up
corresponding to $E_j$.

From (1.4.3) we deduce

\ms\nin
{\bf 1.5.~Proposition.} {\it For $\a={i\over d}+\ell$ with
$i=d-k\in[1,d]$ and $\ell=n-1-p\in[0,n-1]$
$$n_{f,\a}=(-1)^{\ell}\chi(\Y,\Omega_{\Y}^{p}
(\log\Z)\otimes_{\O}\O_{\Y}(-k\H+\msum_{j\in J}\,\lf km_j/d\rf E_j)).
\leqno(1.5.1)$$
Here $\lf\beta\rf:=\max\{k\in\ZZ\mid k\le\beta\}$.}

\ms\nin
{\bf 1.6.~Compatibility with the usual definition.}
The above definition of the spectrum coincides with the definition
using the mixed Hodge structure obtained by $H^ki_0^*\psi_f\Q_X$
where $\psi_f$ denotes the nearby cycle functor (see [11]) and
$i_0:\{0\}\to X_0:=f^{-1}(0)$ denotes the inclusion.
Indeed, by the compatibility of $\psi$ with the direct image
under the blow-up $\rho:\X\to X$ at $0$ (see [26], 2.14), we have
isomorphisms of mixed Hodge modules
$${}^p\!R^k\rho_*\psi_{\f}(\Q_{\X}[n])=\cases{\psi_f(\Q_X[n])&if
$\,k=0$,\cr \,0&if $\,k\ne 0$,}$$
where $\psi$ is shifted by $-1$ so that it preserves the perverse
sheaves.
So we get an isomorphism in the derived category of (algebraic)
mixed Hodge modules
$$\psi_f\Q_X=\R\rho_*\psi_{\f}\Q_{\X}.$$
Note that $Y$ is identified with the exceptional divisor of the
blow-up.
Let $i_Y:Y\to\X_0:=\f^{-1}(0)$ and $a_Y:Y\to pt$ denote the natural
morphisms.
Using the base change of $\rho$ by $i_0$, we then get
$$H^ki_0^*\psi_f\Q_X=H^ki_0^*\R\rho_*\psi_{\f}\Q_{\X}=
H^k(a_Y)_*i_Y^*\psi_{\f}\Q_{\X}.$$
Moreover, if $j':U\to\X_0$ denotes the inclusion, we have by
[6], 4.2
$$(i_Y)_*i_Y^*\psi_{\f}\Q_{\X}=
\R j'_*j'{}^*\psi_{\f}\Q_{\X}.$$
We get thus
$$H^ki_0^*\psi_f\Q_X=H^k(a_U)_*j'{}^*\psi_{\f}\Q_{\X}.$$
So the desired compatibility is reduced to the isomorphism between
the $\l$-eigenspace of the local system $\psi_{\f}\C_{\X}|_U$
and the local system corresponding to $\S^k|_U$ where
$\l=\exp(2\pi\sqrt{-1}k/d)$ since the Hodge filtration is given by
the filtration $\sigma_{\ge p}$ on the logarithmic de Rham complex.
To show the isomorphism of local systems,
it is enough to show the coincidence of the
local monodromies of the two local systems along any irreducible
components
(indeed, this implies the triviality of the tensor product of
one local system with the inverse of the other since a local
system of rank 1 on $U\subset\P^{n-1}$ is trivial if the local
monodromies are).
By the calculation of the residue of the connection in (1.4),
the local monodromy of $\S^k|_U$ along $E_j$ is the multiplication
by $\exp(-2\pi\sqrt{-1}km_j/d)$ where $m_j$ is the multiplicity of
$\f$ along $E_j$.
For $\l=\exp(2\pi\sqrt{-1}k/d)$, it is well-known that the
monodromy of the $\l$-eigenspace of the local system
$\psi_{\f}\C_{\X}|_U$ along $E_j$ is the multiplication by
$\l^{-m_j}$ (since the Milnor fiber is locally defined by
$x_0^dx_j^{m_j}=t$ on a neighborhood of a general point of $E_j$
where $x_0$ defines the exceptional divisor of the blow-up along
$0\in X$), see also [26], 3.3.
So the assertion follows.

\ms\nin
{\bf 1.7.~Weight filtration on the cohomology of the complement.}
From now on, we assume that $D$ is an essential central hyperplane
arrangement.
Let $U=Y\setminus\P(D)$ with $j_U:U\to Y=\P^{n-1}$ the inclusion.
Since $U$ is affine, it is known that
$$(j_U)_!\Q_U[n-1],\q\R(j_U)_*\Q_U[n-1]$$ are perverse sheaves
(see [1]) and underlie mixed Hodge modules (see [26]).
So they have the canonical weight filtration $W$ as perverse
sheaves.
Set
$$\SS(D)=\S(D)\cup\{X\},\q\SS(D)^{(j)}=\{V\in\SS(D)\mid\g(V)=j\},$$
where $\g(V)=\codim\,V$.
Then we can show
$$\eqalign{\Gr^W_{n-1-j}((j_U)_!\Q_U[n-1])&=\mopl_{V\in\SS(D)^{(j)}}
L_V[n-1-j],\cr
\Gr^W_{n-1+j}(\R(j_U)_*\Q_U[n-1])&=\mopl_{V\in\SS(D)^{(j)}}
L_V(-j)[n-1-j],}\leqno(1.7.1)$$
where the $L_V$ are polarized constant variations of Hodge
structures of type $(0,0)$ on $\P(V)\subset\P^{n-1}$.
It is enough to show the first isomorphism since the $L_V$ are
self-dual by the polarization so that the second follows from the
first by duality.

By increasing induction on $j\ge 0$, we define $\Q$-vector spaces
$L_V\,\bigl(V\in\SS(D)^{(j)}\bigr)$ together with morphisms
$$\rho_{V,V'}:L_{V'}\to L_V\q\h{for}\,\,\,V\in\SS(D)^{(j)},\,\,
V'\in\SS(D)^{(j-1)},$$
such that $\rho_{V,V'}=0$ unless $V\subset V'$ as follows:
Set $L_V=\Q$ if $j=0$ or $1$, and $\rho_{V,V'}=id$ if
$\g(V)=1$, $\g(V')=0$.
Assume $L_V$ and $\rho_{V,V'}$ are defined for $\g(V)<j$.
For $V\in\SS(D)^{(j)}$, define
$$\eqalign{&L_V=\Coker\bigl(\msum\,\rho_{V',V''}:
\mopl_{V''\in\SS(D)^{j-2}_V}\,L_{V''}\to
\mopl_{V'\in\SS(D)^{j-1}_V}\,L_{V'}\bigr),\cr
&\q\h{where}\q\SS(D)^{(j)}_{V'}:=\{V\in\SS(D)\mid V\supset V'\}.}
\leqno(1.7.2)$$
The morphism $\rho_{V,V'}:L_{V'}\to L_V$ for $V'\in\SS(D)^{(j-1)}_V$
is given by the composition
$$L_{V'}\to\mopl_{V'\in\SS(D)^{j-1}_V}\,L_{V'}\to L_V,$$
where the first morphism is the natural inclusion, and the second
is the projection to the quotient.
Note that the $L_V$ are Hodge structures of type $(0,0)$,
and they have a canonical polarization using the semi-simplicity
induced by the polarization inductively.

From now on, $L_V$ will be identified with a constant sheaf on
$\P(V)\subset Y=\P^{n-1}$ with stalk $L_V$.
Then $\rho_{V',V}:L_{V'}\to L_V$ is viewed as a morphism of sheaves.
Define a sheaf on $Y$ by
$$\K_Y^j=\mopl_{V\in\SS(D)^{(j)}}\,L_V.$$
We have the morphism $d:\K_Y^{j-1}\to\K_Y^j$ induced by the
$\rho_{V',V}$, and $d\scirc d=0$ by the above construction.
Note that the $\K_Y^j[n-1-j]$ and hence $\K_Y^{\ssb}[n-1]$ are
perverse sheaves on $Y$, see [1].
Then (1.7.1) is reduced to the following lemma since
the weight filtration $W$ is induced by the truncation
$\sigma_{\ge k}$ up to a shift.

\ms\nin
{\bf 1.8.~Lemma.} {\it There is a quasi-isomorphism
$$(j_U)_!\Q_U\simto\K_Y^{\ssb},\leqno(1.8.1)$$
induced by the canonical morphism $(j_U)_!\Q_U\to\Q_Y=\K_Y^0$.}

\ms\nin
{\it Proof.}
Set
$$Y^{(j)}=\mcup_{V\in\SS(D)^{(j)}}\P(V),\q
U^{(j)}=Y^{(j)}\setminus Y^{(j+1)}.$$
Let $\sigma_{\le j}\K_Y^{\ssb}$ denote
the quotient complex of $\K_Y^{\ssb}$ as in [10], 1.4.7, i.e.\
$(\sigma_{\le j}\K_Y)^i=\K_Y^i$ if $i\le j$, and $0$ otherwise.
By increasing induction on $j$ we show
\ms\nin
$(A_j)\q C((j_U)_!\Q_U\to\sigma_{\le j}\K_Y^{\ssb})
[n-2]$ is a perverse sheaf supported on $Y^{(j+1)}$.
\ms\nin
This is clear if $j=0$.
Assume $(A_{j-1})$ holds with $j>0$.
Let $y$ be a general point of $\P(V)$ with $V\in\SS(D)^{(j)}$.
Then $(A_{j-1})$ implies
$$H^k(\sigma_{<j}\K_{Y,y}^{\ssb})=0\q\h{for}\,\,\,k\ne j-1.$$
(Indeed, the restriction of any perverse sheaf to a sufficiently
small Zariski-open subvariety of its support is a local system
shifted by the dimension of the variety.)
Moreover, by (1.7.2), we have the isomorphism as vector spaces
$$H^{j-1}(\sigma_{<j}\K_{Y,y}^{\ssb})=L_V,$$
and this implies the acyclicity of $\sigma_{\le j}\K_{Y,y}^{\ssb}$.
Since the restriction of the cohomology sheaves
${\cal H}^i\sigma_{\le j}\K_Y^{\ssb}$ to $U^{(j)}$ are locally
constant, we see that $\sigma_{\le j}\K_Y^{\ssb}$ is acyclic on
$U^{(j)}$ and hence on $Y^{(1)}\setminus Y^{(j+1)}$ using
$(A_{j-1})$ on $Y^{(1)}\setminus Y^{(j)}$
(since $(j_U)_!\Q_U|_{Y^{(1)}}=0$).
So the shifted mapping cone in $(A_j)$ is supported on $Y^{(j+1)}$,
and it remains to show that the shifted mapping cone is a perverse
sheaf, i.e.\ in the abelian category of perverse sheaves (see [1])
we have
$$\Coker((j_U)_!\Q_U[n-1]\to(\sigma_{\le j}\K_Y^{\ssb})[n-1])=0.$$
But this is clear since the $\Q_{\P(V)}[\dim\P(V)]$ are simple
perverse sheaves so that there are no nontrivial subquotients of
the perverse sheaf $\K_Y^i[n-1-i]$ supported on $Y^{(j+1)}$ if
$i\le j$.
Thus we get $(A_j)$, and (1.8.1) follows by induction.

\ms\nin
{\bf 1.9.~Remark.}
Set $r_V=\rank\,L_V$, and
$\SS(D)_V:=\{V'\in\SS(D)\mid V'\supset V\}$.
By Lemma~(1.8) we have
$$\msum_{V'\in\SS(D)_V}(-1)^{\g(V')}r_{V'}=0.\leqno(1.9.1)$$
This implies that $(-1)^{\g(V)}r_V$ coincides with the
M\"obius function defined by increasing induction on $\g(V)$,
see [24].

Take any $D_k$, and set
$$\SS(D)^{(j)}_{\lkr}=\{V\in\SS(D)^{(j)}\mid V
\not\subset D_k\},\q U_k=Y\setminus D_k=\C^{n-1},$$
with the inclusion $j_k:U\to U_k$.
Then
$$\Gr^W_{n-1+j}(\R(j_k)_*\Q_U[n-1])=\mopl_{V\in\SS(D)_{\lkr}^{(j)}}
L_V|_{U_k}(-j)[n-1-j].$$
The associated spectral sequence degenerates at $E_1$,
since the $U_k\cap\P(V)$ are affine spaces.
This implies that $H^j(U,\Q)$ has type $(j,j)$, and we get
$F^jH^j(U,\C)=P^jH^j(U,\C)=H^j(U,\C)$ for any $j$,
where $P$ is the pole order filtration.
This gives examples where $P\ne F$ locally but $P=F$ globally,
see [13].
The above assertion is compatible with a result of E.~Brieskorn~[2]
that $H^j(U,\Q)$ is generated by logarithmic forms
$$\h{${dg_{i_1}\over g_{i_1}}\wedge\cdots\wedge{dg_{i_j}\over
g_{i_j}}$,}$$
where the $g_i$ are linear functions with constant terms defining
$\P(D_i)\setminus\P(D_k)\subset\C^{n-1}$.
The $E_1$-degeneration also implies a formula in [24]
$$b_k(U)=\msum_{V\in\SS(D)_{\lkr}^{(j)}}\,r_V.\leqno(1.9.2)$$
It is well-known that the Betti numbers $b_k(U)$ are combinatorial
invariants of a hyperplane arrangement, see [24].
We have a refinement as follows.

\ms\nin
{\bf 1.10.~Proposition.} {\it
Set $\S=\S(D)^{\nnc}$, $\S^{D_i}=\{V\in\S\mid V\subset D_i\}$ in
the notation of $(1.1)$.
Then the $b_k(U)$ are determined by
$\S,{\subset},\mu_{\red},\g$ together with
$\S^{D_i}\,(i\in\L)$.}

\ms\nin
{\it Proof.}
We first show that the $r_V$ are determined by the above
combinatorial data.
We have $r_V=1$ for $\SS(D)\setminus\S$ since $\K_{Y,y}^{\ssb}$
for $y\notin\P(D^{\nnc})$ is the standard Koszul complex.
Since the $r_V$ for $V\in\S$ is determined by induction on $\g(V)$
using (1.9.1), it is enough to express
$$\big|\SS(D)_V^{(j)}\setminus\S\big|,$$
using only the combinatorial data as above.
Set $I(V)=\{i\in\L\mid D_i\supset V\}$.
Note that $I(V)$ is determined by the $\S^{D_i}\,(i\in\L)$ if
$V\in\S$. Set
$$\eqalign{S(\L)^{(j)}&=\{I\subset\L\mid|I|=j\},\q
S(\L)_V=\{I\subset\L\mid I\subset I(V)\},\cr
S^{\nnc}(\L)&=\{I\subset\L\mid I=I(V')\,\,\,\h{for some}\,\,
V'\in\S\}.}$$
Then we have the identification
$$\SS(D)_V^{(j)}\setminus\S=S(\L)^{(j)}\cap S(\L)_V\setminus
S^{\nnc}(\L),$$
and the assertion follows.
Thus the $r_V$ are calculated by induction on $\g(V)$ using
only the above combinatorial data.

Since the Betti number is calculated by using (1.9.2), it is enough
to express
$$\big|\SS(D)_{\lkr}^{(j)}\setminus\S\big|,$$
by using only the combinatorial data as above.
So the assertion follows since
$$\SS(D)_{\lkr}^{(j)}\setminus\S=S(\L)^{(j)}\cap S(\L)_{\lkr}
\setminus S^{\nnc}(\L),$$
where $S(\L)_{\lkr}=\{I\subset\L\mid k\notin I\}$.
This finishes the proof of Proposition~(1.10).

\bs\bs
\centerline{\bf 2. Canonical embedded resolution}

\bs\nin
The material in this section is treated in a much more general
situation in [8].
For the convenience of the reader we treat it under the assumption
that $\S$ is stable by intersection.
This hypothesis is satisfied in our case, and simplifies certain
arguments very much.

\ms\nin
{\bf 2.1.~Construction.}
Let $\S$ be a finite set of proper vector subspaces of the vector
space $X=\C^n$ which is stable by intersection
(i.e. $V\cap V'\in\S$ if $V,V'\in\S$) and such that $0\in\S$.
We have a function $\g:\S\to\N$ associating the codimension of
$V$.
There is a natural order on $\S$ defined by the inclusion relation.
We say that $\S$ and $\S'$ are combinatorially equivalent if there
is a one-to-one correspondence $\S\to\S'$ as ordered sets in a
compatible way with $\g$.
Note that a {\it nested} subset of $\S$ in the sense of [8] is
always linearly ordered by the inclusion relation in our paper
since $\S$ is stable by intersection.

Let $Y=\P^{n-1}$.
For a vector subspace $V\subset X=\C^{n}$, its corresponding
subspace of $Y$ will be denoted by $\P(V)$.
For $\S$ as above, there is a sequence of blowing-ups
$$\rho_i:Y_{i+1}\to Y_i\q\h{for}\,\,\, 0\le i<n-2,$$
whose center $C_i$ is the {\it disjoint} union of the proper
transforms of $\P(V)$ for $V\in\S$ with $\dim\P(V)=i$,
where $Y_0=Y$.
Note that we cannot restrict to the dense edges as in [30] because
this is not adequate for our inductive argument.

Set $\Y=Y_{n-2}$ with $\rho:\Y\to Y$ the composition of the
$\rho_i$.
We will sometimes denote $\Y$ by $Y^{\S}$.
This applies to $\P(V)^{\S^V}$ where $Y$ and $\S$ are
replaced by $\P(V)$ and $\S^V$ respectively.
Here we define for $V\in\S$
$$\S^V=\{V'\in\S\mid V'\ssn V\},\q
\S_V=\{V'\in\S\mid V'\supset V\}.$$
For $V\in\S\setminus\{0\}$ with $\dim\P(V)=i$, let $C_{V,j}$
denote the proper transform of $C_{V,0}:=\P(V)\subset Y_0$ in
$Y_j$ for $1\le j\le i$.
Let $E_{V,i+1}$ be the exceptional divisor of the blow-up
along $C_{V,i}$ which is an irreducible component of $C_i$
and is identified with $\P(V)^{\S^V}$ (which is defined above).
Let $E_{V,j}$ be the proper transform of $E_{V,i+1}$ in $Y_j$ for
$i+1< j\le n-2$.
Finally, set $E_V=E_{V,n-2}$ if $\dim\P(V)<n-2$, and $E_V=C_{V,n-2}$
if $\dim\P(V)=n-2$.
For $V=0$, let $E_0$ denote the proper transform $\H$ in $\Y$ of a
{\it general\ } hyperplane $H$ of $\P^{n-1}$.
(It is known that the divisor class group of $\Y$ is generated by
$E_V$ for $V\in\S$ with $\dim\P(V)<n-2$.)

\ms\nin
{\bf 2.2.~Remarks.}
(i) For smooth varieties $X\supset Y\supset Z$ in general,
the proper transform of $Y$ by the blow-up of $X$ along $Z$ is the
blow-up of $Y$ along $Z$.

\ms
(ii) For any linear subspaces $L,L'$ of affine space such that
$L\cap L'\ne L,L'$, the proper transforms of $L$ and $L'$ by the
blow-up along $L\cap L'$ do not intersect.

\ms\nin
{\bf 2.3.~Proposition.}
{\it The union of $E_V$ for $V\in\S\setminus\{0\}$ is a divisor
with normal crossings on $\Y$, and the intersection of $E_{V_k}$
for $V_k\in\S\setminus\{0\}$ with $1\le k\le r$ is empty unless
$V_1\subset\cdots\subset V_k$ up to a permutation.
}

\ms\nin
{\it Proof.}
The last assertion follows from Remark~(2.2)(i) because $\S$ is
stable by intersection.
For the first assertion we take local coordinates
$x_1,\dots,x_{n-1}$ such that $V_k$ is given by $x_i=0$ for
$i>d_k$ where $d_k=\dim V_k$.
Then local coordinates $y_1,\dots,y_{n-1}$ of the blow-up along
$V_1$ are given by $y_i=x_i$ if $i\le d_1$ or $i=i_1$, and by
$y_i=x_i/x_{i_1}$ otherwise. Here the exceptional divisor of the
blow-up is given by $y_{i_1}=0$ and $i_1$ is an integer in
$(d_1,d_2]$, because we have to consider the proper transforms of
the $V_k$ for $k\ge 2$, which are given by $y_i=0$ for $i>d_k$
if $i_1\in(d_1,d_2]$.
So the assertion follows by repeating this construction.

\ms\nin
{\bf 2.4.~Proposition.}
{\it The center $C_i$ of the blow-up $\rho_i$ is the disjoint union
of $C_{V,i}=\P(V)^{\S^V}$ for $V\in\S$ with $\dim\P(V)=i$, and
$C_{V,j+1}\,(1\le j<i)$ is the blow-up of $C_{V,j}$ along the
disjoint union of $C_{V',j}$ for $V'\in\S^V$ with $\dim\P(V')=j$.
}

\ms\nin
{\it Proof.}
This follows from Remarks~(2.2)(i) and (ii) using the calculation
in (2.3).

\ms\nin
{\bf 2.5.~Proposition.}
{\it Let $V,V'\in\S$ such that $\dim\P(V)=i<\dim\P(V')=i'$.
Then $C_{V',i'}$ is not contained in $E_{V,i'}$, and hence
$E_{V,j}$ inductively coincides with the total transform of
$E_{V,i}$ for $j>i$.}

\ms\nin
{\it Proof.}
This follows from the above arguments, because $\S$ is stable by
intersections and the blow-ups are done by increasing induction on
the dimension of the center.
More precisely, we have $C_{V',i'}\cap E_{V,i'}=\emptyset$
in the case $V\not\subset V'$.
In the other case, repeating the above construction with $Y$
replaced by $\P(V')$, define $C'_{V,i}$, $E'_{V,i+1}$, and
$E'_{V}$ associated to $\P(V)\subset\P(V')$ as in $(2.1)$ above,
i.e.\ $E'_{V,i+1}$ is the exceptional divisor of the blow-up
along the center $C'_{V,i}$, and $E'_{V}$ is the proper
transform of $E'_{V,i+1}$ in $\P(V')^{\S^{V'}}=C_{V',i'}$, where
the upper script $'$ of $C'_{V,i}$, $E'_{V,i+1}$ and $E'_{V}$
means that the construction is done for $\P(V')$ instead of $Y$.
Then $C_{V',i'}\cap E_{V,i'}=E'_{V}$.
This finishes the proof of Proposition~(2.5).

\ms\nin
{\bf 2.6.~Proposition.}
{\it For $V\in\S$, $E_V$ depends only on $\S^V$ and $\S_V$, and
there is a canonical decomposition
$$E_V=\P(V)^{\S^V}\times\P(X/V)^{\S_V},$$
where $\P(V)^{\S^V}$, $\P(X/V)^{\S_V}$ are the successive blow-ups
of $\P(V)$ and $\P(X/V)$ respectively associated with $\S^V$ and
$\S_V$ as in {\rm (2.1)}.}

\ms\nin
{\it Proof.}
The first assertion is clear because $\S$ is stable by intersections.
Let $r$ be the codimension of $\P(V)$ in $Y$ (i.e. $r=n-1-i$).
Taking $r$ general hyperplanes containing $V$, and considering their
proper transforms whose intersection is $C_{V,i}$, we see that the
tensor of the conormal bundle of $C_{V,i}$ with some line bundle is
a trivial vector bundle.
Hence $E_{V,i+1}$ is a trivial $\P^{r-1}$-bundle over $C_{V,i}$
with a canonical projection to $\P^{r-1}$ (which is independent of
the choice of the hyperplanes up to the action of $PGL(r-1,\C)$).
So the assertion follows, because $\S_V$ is identified with a set
of vector subspaces of $X/V$.

\ms\nin
{\bf 2.7.~Proposition.}
{\it Let $V_k\in\S$ for $1\le k\le r$ such that $0\ne V_1
\ssn\cdots\ssn V_r$. Then
$$\mcap_{1\le k\le r}\,E_{V_k}=\mprod_{0\le k\le r}\,
\P(V_{k+1}/V_k)^{\S_k},$$ where
$V_0=0,V_{r+1}=X$, and $\S_k=\{V/V_k\mid V\in\S\,\,\,
\h{with}\,\,\,V_k\subset V\ssn V_{k+1}\}$.
}

\ms\nin
{\it Proof.}
Let $V=V_r$ and $i=\dim\P(V)$. Then $E_{V}\cap E_{V_k}$
is the pull-back of $C_{V,i}\cap E_{V_k,i}$ by the projection
$E_{V}\to C_{V,i}$ for $1\le k\le r-1$.
So the assertion follows from (2.6) by induction on $r$,
where the inductive hypothesis is applied to $C_{V,i}
=\P(V)^{\S^V}$ and $C_{V,i}\cap E_{V_k,i}\,(1\le k\le r-1)$
which are calculated in the proof of Proposition~(2.5).

\ms\nin
{\bf 2.8.~Proposition.}
{\it For $V'\in\S^V$ {\rm (}resp. $\S_V)$ such that $V'\ne V$,
the pull-back of $E_{V'}$ to $E_V$ coincides with the pull-back
of $E'_{V'}$ on $\P(V)^{\S^V}$ by $pr_1$ {\rm (}resp. of $E''_{V'}$
on $\P(X/V)^{\S_V}$ by $pr_2)$.
Here $\P(V)^{\S^V}$ and $\P(X/V)^{\S_V}$ are as in {\rm (2.6)}, and
$pr_i$ denote the projection to the $i$-th factor.
For $V'=V$, the pull-back of $E_{V}$ to $E_V$ as a divisor class
is given by
$$pr_1^*(\H'-\msum_{V'\in\S^V}\,E'_{V'})-pr_2^*\H'',$$
where $\H'$ {\rm (}resp. $\H'')$ is the proper transform of a
general hyperplane of $\P(V)$ {\rm (}resp. $\P(X/V))$, and it is
zero if $\dim\P(V)=0$ {\rm (}resp. $\dim\P(X/V)=0)$.
}

\ms\nin
{\it Proof.}
The first assertion follows from (2.5--6).
For the second, take a general hyperplane of $Y$
intersecting $\P(V)$ transversally and also a hyperplane of $Y$
containing $\P(V)$ and corresponding to a general hyperplane of
$\P(X/V)$.
Then the assertion follows by considering the difference between
their pull-backs to $\Y$.

\bs\bs
\centerline{\bf 3. Proofs of Theorems~3 and 5}

\bs\nin
{\bf 3.1.~Proof of Theorem~3.}
This follows from Proposition~(1.5) together with the Riemann-Roch
theorem for surfaces (see e.g.\ [17], Example~15.2.2)
$$\h{$\chi(\O_{\Y}(D_k))={1\over 2}\,D_k\cdot(D_k-K_{\Y})+
\chi(\O_{\Y})$}.\leqno(3.1.1)$$
Here $\Y$ is the blow-up of $Y=\P^2$ along the points of
$\P(D^{\nnc})$ in the notation of (0.1), and $D_k$ is a divisor on
it (which will be defined later).

By (1.4.1) we may assume
$$\a={i\over d}+\ell\q\h{with}\,\,\,i=d-k\in[1,d],\,\,\ell=2-p
\in[0,2],$$
since $n_{f,\a}=0$ for the other $\a$.
Consider first the case $\ell=0$. We have
$$\Omega_{\Y}^2(\log\Z)=\O_{\Y}((d-3)\H+\msum_{j\in J''}\,
(2-m_j)E_j)),$$
since $\Omega_{\Y}^2=\rho^*\Omega_Y^2\otimes_{\O}\O_{\Y}
(\msum_{j\in J''}\,E_j)$ and $\Z_{\red}=\rho^*Z+\msum_{j\in J''}\,
(1-m_j)E_j$ in the notation of (1.4).
So we apply the Riemann-Roch theorem (3.1.1) to
$$D_k=(d-k-3)\H+\msum_{j\in J''}(-m_j+\lf km_j/d\rf+2)E_j,$$
using Proposition~(1.5) for $p=2$.
We have $D_k=\msum_j\,A_jE_j+C\H$ with
$$A_j=2+\lf-im_j/d\rf=2-\lc im_j/d\rc,\q C=i-3,$$
and $K_{\Y}=-3\H+\msum_j\,E_j$. So we get
$$D_k^2=-\msum_j\,A_j^2+C^2,\q D_k\cdot K_{\Y}=-\msum_j\,A_j-3C,$$
since $\H^2=1$, $E_j^2=-1$, and $\H,E_j$ are orthogonal to each
other.
These imply the first equality by (3.1.1) since $\chi(\O_{\Y})=1$.

The argument is similar for the last equality where $\ell=2$,
$p=0$, $D_k=\msum_j\,A_jE_j+C\H$ with
$$A_j=\lf km_j/d\rf=m_j-\lc im_j/d\rc,\q C=-k=i-d.$$
Note that the reduced cohomology is used for the definition
of spectrum, and the difference corresponds to $\delta_{i,d}$
in the case $i=d$.

By the identity ${a+b\choose 2}-{a\choose 2}-{b\choose 2}=ab$,
the middle equality for $\ell=1$ follows from the others since
we have by (1.4.2)
$$\h{$\sum_{\ell=0}^2n_{f,{i\over d}+\ell}=\chi(U)-\delta_{i,d}\q
\h{with}\q\chi(U)={d-2\choose 2}-\sum_{m\ge 3}\nu^{(2)}_m
{m-1\choose 2}$}.$$
Here the first equality is clear by the definition of spectrum.
The last equality is shown by using a small deformation
to a generic central arrangement $D'$ where $\P(D')$ is a divisor
with normal crossing so that
$$\h{$\chi(\P^2\setminus\P(D'))=3-2d-{d\choose 2}={d-2\choose 2}$}.$$
The difference of the local Euler characteristics of $\P(D)$ and
$\P(D')$ at each point of $\P(D)$ with multiplicity $m$ is given by
$1-\bigl(m-{m\choose 2}\bigr)={m-1\choose 2}$.
So the assertion follows.

\ms\nin
{\bf 3.2.~Generic case.}
Assume $D$ is a generic central hyperplane arrangement, i.e.\
$\P(D)\subset\P^{n-1}$ is a divisor with normal crossings.
In this case it is known ([27], Cor.~1) that
$$\h{$n_{f,\a}=n_{f,n-\a}={i-1\choose n-1}\,\,\,$ for
$\,\a=i/d<1$},$$
where $n_{f,\a}=0$ for $d\a\notin\ZZ$.
For $\a=i\in\ZZ$, we have by 5.6.1 in loc.~cit.
$$\h{$n_{f,i}=(-1)^{i-1}{d-1\choose n-i}\,\,\,$ for
$\,1\le i\le n-1$}.$$
It is possible to calculate $n_{f,\a}$ for any $\a$ using
Proposition~(1.5) and the Bott vanishing theorem.

\ms\nin
{\bf 3.3.~{\Mustata}'s formula.}
In the notation of (1.1), set $\S'=\S(D^{\nrnc})$.
For each $V\in\S'$, let $\I_V\subset\C[X]$ be the reduced ideal of
$V$.
{\Mustata}'s formula [22] states that for any $\a>0$
$$\J(X,\a D)=\mcap_{V\in\S'}\,\I_V^{\lf\a\mu(V)\rf-
\g(V)+1}.\leqno(3.3.1)$$
In the nonreduced case this is due to Z.~Teitler [32] (see also
[27], 2.2).

If $V=0$, then we have for $\a=j/d$ with $j\in[1,d]\cap\ZZ$
$$\I_0^{\lf(\a-\varepsilon)\mu(0)\rf-(\g(0))+1}=
\I_0^{j-n},$$
for $0<\varepsilon\ll 1/d$, since $\g(0)=n$ and
$\mu(0)=\deg D=d$.
So we get (0.1) where $n=3$, see also [22], Cor.~2.1.

\ms\nin
{\bf 3.4.~Proof of Theorem~5.}
By (1.4.3--5) for $p=n-1$, we get for $\a={i\over d}\in(0,1]$
with $i=d-k\in[1,d]$
$$\eqalign{n_{f,\a}&=\dim\Gamma\bigl(\Y,\O_{\Y}
((i-n)\H+\msum_{j\in J}\,(c_j-\lc im_j/d\rc)E_j)\bigr)\cr
&=\dim\bigl(\mcap_{V\in\S'\setminus\{0\}}\,\I_V^{\lceil\a\mu(V)
\rceil-\g(V)}\cap\C[X]_{i-n}\bigr),}\leqno(3.4.1)$$
where the second equality is shown by using the injection
$$\eqalign{&\Gamma\bigl(\Y,\O_{\Y}((i-n)\H+\msum_j\,
\bigl(c_j-\lc im_j/d\rc\bigr)E_j\bigr)\cr
&\subset\Gamma(\Y,\O_{\Y}((i-n)\H))=\Gamma(Y,\O_Y(i-n))=
\C[X]_{i-n}.}$$
Indeed, for $g\in\C[X]_{i-n}=\Gamma(\Y,\O_{\Y}((i-n)\H))$, the
condition $g\in\I_V^k$ corresponds to that $\pi^*g\in\I_{E_j}^k$
if $E_j$ corresponds to $V$, where $\I_{E_j}$ is the ideal of $E_j$
and
$$k=\lc im_j/d\rc-c_j=\lc\a\mu(V)\rc-\g(V).$$
Note that we may have $c_j-\lc im_j/d\rc>0$ only in the case
$c_j\ge 2$ so that $g$ cannot have a pole.
Thus the assertion follows.

\ms\nin
{\bf 3.5.~Relation with $b$-functions.}
It does not seem easy to get an explicit formula for the
jumping coefficients and the spectrum of a hyperplane arrangement
in the general case.
However, it seems to be more difficult to calculate the roots of
the $b$-function of a hyperplane arrangement except for the case
of a generic central arrangement [33], see also [27].
The relation between the jumping coefficients $\JC(D)$ and the
roots of the $b$-function $R_f$ is quite complicated although
there is an inclusion relation
$$\JC(D)\cap(0,1)\subset R_f\cap(0,1),$$ as is shown in [14].
The converse inclusion holds under some conditions, see [27].
However, it does not hold without these conditions
as is shown by the following.

\ms\nin
{\bf 3.6.~Example.} Let $n=3$, $d=7$ and
$f=(x^2-y^2)(x^2-z^2)(y^2-z^2)z$.
By (0.1) or Theorem~3 we have ${5\over 7}\notin\JC(D)$, but
${5\over 7}\in R_f$, see [28].
In this case
$$\dim\C[x,y,z]_2=\#\P(D^{\nnc})=6.\leqno(3.6.1)$$
So we have to prove the non-degeneracy of some matrix to show
the non-existence of a nontrivial homogeneous polynomial of degree
2 vanishing at all the 6 points if we want to show that
${5\over 7}\notin\JC(D)$ using (0.1).
Note that $n_{f,5/7}=0$ by Theorem~3 where $d=7$, $\nu^{(2)}_3=6$
and $\nu^{(2)}_m=0\,(m>3)$.
This implies that ${5\over 7}\notin\JC(D)$ by Proposition~(4.2).

\vfill\eject
\centerline{\bf 4. Proofs of Theorems~1--4 by induction}

\bs\nin
{\bf 4.1.~Isolated jumping coefficients.}
Let $X$ be a smooth variety (or a complex manifold), and $D$ be a
divisor on it.
Let $\J(X,\a D)$ denote the multiplier ideals for $\a>0$, see [21].
The graded pieces of the multiplier ideals are defined for
$\a > 0$ with $0<\varepsilon\ll 1$ by
$$\G(X,\a D)=\J(X,(\a-\varepsilon)D)/\J(X,\a D).$$
The jumping coefficients are the rational numbers $\a$ such that
$\G(X,\a D)\ne 0$.
We say that $\a$ is an {\it isolated\ } jumping coefficient at
$x$ if $\G(X,\a D)$ is supported on $x$.

If $D$ has an isolated singularity at $x$ and is defined locally
by $f$, then the coefficient $n_{f,\a}$ of the spectrum
$\Sp(f)=\sum_{\a>0}n_{f,\a}t^{\a}$ for $\a\in(0,1)$ is given
(see [3]) by
$$n_{f,\a}=\dim\G(X,\a D)_x.\leqno(4.1.1)$$

\ms\nin
{\bf 4.2.~Proposition.} {\it
The assertion $(4.1.1)$ holds by assuming only that $\a\in(0,1)$ is
an isolated jumping coefficient at $x$.}

\ms\nin
{\it Proof.}
By [6], we have a canonical isomorphism
$$\G(X,\a D)=F_{-n}\Gr_V^{\a}\B_f,$$
where $\Gr_V^{\a}\B_f$ coincides with the $\l$-eigenspace of
$\psi_f\O_X$ by the action of the monodromy where $\l=e^{-2\pi i\a}$.
Here we have to show that $F_{-n}$ does not change by taking the
pull-back by $i_x:\{x\}\to X$ as in (1.6).
Choosing local coordinates $x_1,\dots,x_n$ of $(X,x)$ and using
[26], 2.24, the pull-back $i_x^*$ is given by iterating
$$i_k^*=C({\rm can}:\psi_{x_k,1}\to\varphi_{x_k,1}).$$
For the underlying filtered left $D$-modules $(M,F)$, the last
functor is given by the mapping cone of
$$\d_{x_k}:\Gr_{V_k}^1(M,F[1])\to\Gr_{V_k}^0(M,F),$$
where $V_k$ is the $V$-filtration along $x_k=0$ and
$x_k\d_k-\a$ is nilpotent on $\Gr_{V_k}^{\a}$.
Since ${\rm supp}\,F_{-n}M=\{x\}$, we see that $F_{-n}M$ is
contained in
$$(i_x)_*H^0i_x^!M=\Gamma_{[x]}M\subset M,$$
which underlies a mixed Hodge module, and is isomorphic to
$\mopl_i\,(\C[\d_1,\dots,\d_n],F[p_i])$ with $p_i\in\ZZ$.
(Indeed, it is the direct image as a $D$-module of the filtered
vector space $H^0i_x^!(M,F)$ by the closed embedding $i_x$.
Note that $\C[\d_1,\dots,\d_n]$ has the Hodge filtration $F$ by the
degree of polynomials in $\d_i$.)
So $F_{-n}M$ does not change by passing to $\Gr_{V_k}^0(M,F)$
inductively.
Since $F_{-n}=0$ on $(M,F[1])$, this implies the desired result.

\ms\nin
{\bf 4.3.~Proof of Proposition~1.}
Assume $\a\in\JC(D)\cap(0,1)$.
It is well-known that $\G(X,\a D)=0$ for $\a\in(0,1)$ if $D$ is a
reduced divisor with normal crossings.
So the support of $\G(X,\a D)$ is a union of $V\in\S(D)^{\nrnc}$,
since $D$ is locally trivial along a non-empty Zariski open subset
of $V$.
Restricting $D$ to an affine subspace which is transversal to
the Zariski open subset of $V\in\S(D)^{\nrnc}$ and has
complementary dimension with $V$, it is enough to consider
the case of isolated jumping coefficients, since $D$ is locally
the product of its restriction to the transversal space with $V$.
(But this does not mean that the singular points of $D$ are
$0$-dimensional.)
Then we get the non-vanishing of $n_{f_{X/V},\a}$ by
Proposition~(4.2) since $f_{X/V}$ is identified with the defining
polynomial of the restriction of $D$ to the transversal space.
So Proposition~1 follows since the opposite implication is
well-known, see [3].

\ms\nin
{\bf 4.4.~Theorem.} {\it
In the notation of $(1.1)$ and $(2.1)$, set
$\S=\S(D)^{\nnc}$ and
$\S^{D_i}=\{V\in\S\mid V\subset D_i\}$.
Let $\Z$ be as in $(1.4)$.
For $a=(a_V)_{V\in\S}\in\ZZ^{\S}$ define
$$\Phi_{\S}^p(a)=\chi(\Y,\Omega_{\Y}^p(\log\Z)\otimes_{\O}\O_{\Y}
(\msum_{V\in\S}\,a_VE_V)).$$
Then $\Phi_{\S}^p(a)$ is a polynomial in $a=(a_V)_{V\in\S}$ whose
coefficients are rational numbers and are determined by the
combinatorial data of $D$.
More precisely, it depends only on $\S,{\subset},\mu_{\red},\g$
together with $\S^{D_i}\,(i\in\L)$ in the notation of
$(1.1.1)$.
If $p=0$, then $\Phi_{\S}^0(a)$ depends only on the weak
combinatorial data.}

\ms\nin
{\it Proof.}
We show this by increasing induction on $n=\dim X\ge 2$.
First we show the assertion on $\Phi_{\S}^p(a)$ for any $p$.
If $n=2$, the assertion is trivial by the Riemann-Roch theorem
for curves.
Here $\S=\{0\}$ and the number of the points of $\Z$ is enough for
the calculation.

Assume $n>2$, and set
$$M(a)=\Omega_{\Y}^p(\log\Z)\otimes_{\O}\O_{\Y}(\msum_{V\in\S}\,
a_VE_V).$$
Here we may assume $p<\dim \Y$ since the case $p=\dim \Y$ is
reduced to the case $p=0$.
Take some $V\in\S$, and set $E=E_V$ if $V\ne 0$.
In the case $V=0$, set $E=\H$ which is the pull-back of a
sufficiently general hyperplane.
We have a short exact sequence
$$0\to M(a-{\bf 1}_V)\to M(a)\to M(a)\otimes_{\O_{\Y}}\O_{E}\to 0,$$
where ${\bf 1}_V\in\ZZ^{\S}$ is defined by $({\bf 1}_V)_{V'}=0$ for
$V'\ne V$ and $({\bf 1}_V)_V=1$. So we get
$$\Phi_{\S}^p(a)-\Phi_{\S}^p(a-{\bf 1}_V)=
\chi(E,M(a)\otimes_{\O_{\Y}}\O_{E}).\leqno(4.4.1)$$
Using the identity ${x\choose k}-{x-1\choose k}={x-1\choose k-1}$
as polynomials in $x$ where $k\in\ZZ_{>0}$
(see also [18], I, Prop.~7.3(a)),
it is enough to show that (4.4.1) is a polynomial determined by the
combinatorial data.

We consider first the case $V\ne 0$.
Let $N^*_E$ denote the conormal bundle of $E\subset\Y$.
This is the restriction of the line bundle $\O_{\Y}(-E)$,
and the restriction as a divisor is calculated in (2.8).
Let $\Z'$ be the closure of $\Z\setminus E$.
Set $\Z'_E=\Z'\cap E$.
There is a commutative diagram of exact sequences
$$\matrix{&&0&\to&\O_{\Y}(-E)&=&\O_{\Y}(-E)&\to&0\cr
&&\downarrow&&\cap&&\cap\,\,\cr
0&\to&\Omega_{\Y}^1(\log\Z')&\to&M&\to&\O_{\Y}&\to&0\cr
&&||&&\don&&\don\,\,\cr
0&\to&\Omega_{\Y}^1(\log\Z')&\to&\Omega_{\Y}^1(\log\Z)&\to&\O_E&
\to&0}$$
where $M:=\Ker(\Omega_{\Y}^1(\log\Z)\oplus\O_{\Y}\to\O_E)$.
Taking the pull-back by $E\to\Y$, we get short exact sequences
$$\eqalign{0\to N_E^*\to M|_E&\to\Omega_{\Y}^1(\log\Z)|_E\to 0,\cr
0\to\Omega_{\Y}^1(\log\Z')|_E\to M|_E&\to\O_E\to 0.}$$
We have also a short exact sequence
$$0\to N_E^*\to\Omega_{\Y}^1(\log\Z')|_E\to\Omega_E^1(\log(\Z'_E)
\to 0.\leqno(4.4.2)$$
These imply the equalities in the Grothendieck group
$$\eqalign{\lb\mwedge^pM|_E\rb&=\lb\Omega_{\Y}^p(\log\Z)|_E\rb+
\lb N_E^*\otimes\Omega_{\Y}^{p-1}(\log\Z)|_E\rb,\cr
\lb\mwedge^pM|_E\rb&=\lb\Omega_{\Y}^p(\log\Z')|_E\rb+
\lb\Omega_{\Y}^{p-1}(\log\Z')|_E\rb\cr
&=\lb\Omega_E^p(\log\Z'_E)\rb+\lb N_E^*\otimes\Omega_E^{p-1}
(\log\Z'_E)\rb\cr
&\q+\lb\Omega_E^{p-1}(\log\Z'_E)\rb+\lb N_E^*\otimes
\Omega_E^{p-2}(\log\Z'_E)\rb.}$$
So we get by increasing induction on $p$
$$\lb\Omega_{\Y}^p(\log\Z)|_E\rb=\lb\Omega_E^p(\log\Z'_E)\rb+
\lb\Omega_E^{p-1}(\log\Z'_E)\rb.\leqno(4.4.3)$$
The assertion on (4.4.1) is thus reduced to that
$$\chi(E,\Omega_E^p(\log\Z'_E)\otimes_{\O}\O_{\Y}(\msum_{V\in\S}\,
a_VE_V)|_E)$$
is a polynomial determined by the combinatorial data.
By Proposition~(2.6) we have
$$E=\P(V)^{\S^V}\times\P(X/V)^{\S_V},$$
and the restriction of $\O_{\Y}(\msum_{V\in\S}\,a_VE_V))$ to $E$
is calculated by Proposition~(2.8).
Moreover we have the decomposition
$$\Z'_E=pr_1^*Z_1+pr_2^*Z_2,$$
where $Z_1$ is given by $E'_{V'}$ with $V'\subset V$, and $Z_2$
is given by $E''_{D_i/V},E''_{V'/V}$ with $D_i,V'\supset V$.
Here $E''_{D_i/V}\subset\P(X/V)^{\S_V}$ is the proper transform of
$\P(D_i/V)\subset\P(X/V)$, and the associated combinatorial data
$(\S_V)^{D_i/V}$ is given by $\S_V\cap\S^{D_i}$.
So the assertion after (4.4.1) for $V\ne 0$ follows from the
inductive assumption using the K\"unneth-type decomposition of
$\Omega_E^p(\log\Z'_E)$.

For $V=0$ we have a similar assertion since $\H$ intersects
transversally $E_V$ for every $V\in\S\setminus\{0\}$.
Using an exact sequence similar to (4.4.2), we get instead of
(4.4.3)
$$\lb\Omega_{\Y}^p(\log\Z)|_{\H}\rb=\lb\Omega_{\H}^p(\log\Z_{\H})
\rb+\lb N_{\H}^*\otimes\Omega_{\H}^{p-1}(\log\Z_{\H})\rb.
\leqno(4.4.4)$$
Applying the inductive hypothesis, we can then calculate
$$\chi(\H,\Omega_{\Y}^p(\log\Z)|_{\H}\otimes_{\O}\O_{\Y}
(\msum_{V\in\S}\,a_VE_V)|_{\H}),$$
where $\S^{D_i}$ and $\S$ are replaced by those obtained by
deleting the $1$-dimensional $V$ (i.e.\ $\dim\P(V)=0$).
So the assertion after (4.4.1) for $V=0$ is also proved.

Thus the assertion is reduced to the case $a=0$,
and we have to show that
$$\Phi_{\S}^p(0)=\chi(\Y,\Omega_{\Y}^p(\log\Z))$$ depends only on
the combinatorial data.
It is known that each $H^j(U)$ is generated by products of
logarithmic 1-forms (see [2]), and hence has type $(j,j)$.
Then the assertion that the Hodge numbers of $U=\Y\setminus\Z$ depend
only on the combinatorial data is equivalent to a similar assertion
for the Betti numbers.
So the assertion follows from Proposition~(1.10).

The argument is similar and easier for the assertion on
the weak combinatorial data in the case $p=0$,
since we do not have to treat the logarithmic forms by the
isomorphism $\Omega_{\Y}^0(\log\Z)=\O_{\Y}$.
This finishes the proof of Theorem~(4.4).

\ms\nin
{\bf 4.5.~Proofs of Theorems~1--2 by induction.}
Theorem~1 follows from Theorem~(4.4) and Proposition~(1.5).
If $D$ is reduced, then $\lf km_j/d\rf=0$ for $j\in J'$ in the
notation of (1.4), and $c_j-m_j=0$ for $j\in J'$ in (1.4.5).
Using the second equality of (1.4.5) for $p=n-1$,
Theorem~2 then follows from Theorem~(4.4) and Proposition~(1.5).

\ms\nin
{\bf 4.6.~Remark.}
In the original version [29], the argument in the proof of
Theorem~(4.4) treated $\Omega_{\Y}^p$ and not
$\Omega_{\Y}^p(\log\Z)$.
By this method we have to take the graded pieces of the weight
filtration $W$ on the logarithmic forms, and the argument becomes
more complicated.
The above proof of Theorem~1 was obtained after the new proof
in the next section appeared.

We can calculate examples using the method in this section
as is shown below.

\ms\nin
{\bf 4.7.~Proof of Theorem~3 by induction.}
Let $\Phi^p(A,C)$ denote $\Phi_{\S}^p(a)$ in (4.4) where $A=(A_j)$,
and the $a_V$ are denoted by $A_j$ or $C$ depending on whether
$\dim V=1$ or $0$.
Let $E_j$ denote the exceptional divisor corresponding to $A_j$.
Then
$$\Phi(A,C)=\chi(\Y,\E)\q\h{with}\q\E=\O_{\Y}(\msum_jA_jE_j+C\H).$$
We have $\R\pi_*\O_{\Y}=\O_Y$ where $\pi:\Y\to Y$ since $Y$
is nonsingular.
So we first get
$$\h{$\Phi^0(0,C)={C+2\choose 2}$}.$$
Fix $j$, and let ${\bf 1}_j$ be as in the proof of Theorem~(4.4).
We have $\O_{\Y}(E_j)|_{E_j}=\O_{E_j}(-1)$ where $E_j=\P^1$.
(This is shown by using the total transform of a hyperplane passing
through the center of the blow-up.) Hence
$$\Phi^0(A,C)-\Phi^0(A-{\bf 1}_j,C)=\chi(E_j,\O_{E_j}(-A_j))=
1-A_j.$$
Thus we get
$$\h{$\Phi^0(A,C)={C+2\choose 2}-\sum_j{A_j\choose 2}$}.$$
This implies the assertions for $\a\in(0,1]\cup(2,3]$ by
setting $A_j,C$ as in (3.1).

For $p=1$, we have by (4.4.4) applied to a general $H=\P^1\subset Y$
$$\eqalign{\Phi^1(0,C)-\Phi^1(0,C-1)&=
\chi(\P^1,\Omega_{\P^1}^1\otimes\O_{\P^1}(C+d))+
\chi(\P^1,N_H^*\otimes\O_{\P^1}(C))\cr
&=(C+d-1)+C=2C+d-1.}$$
Since $\Phi^1(0,0)=b_1(U)=d-1$, we get
$$\Phi^1(0,C)=C^2+dC+d-1.$$
Fix now $j$.
We have by (4.4.3)
$$\eqalign{\Phi^1(A,C)-\Phi^1(A-{\bf 1}_j,C)&=
\chi(\P^1,\Omega_{\P^1}^1\otimes\O_{\P^1}(m_j-A_j))+
\chi(\P^1,\O_{\P^1}(-A_j))\cr
&=(m_j-A_j-1)+(1-A_j)=m_j-2A_j.}$$
We get thus
$$\Phi^1(A,C)=\msum_j\,(-A_j^2-A_j+m_jA_j)+C^2+dC+d-1.$$
Here $A_j=m_j-\lc im_j/d\rc$, $C=i-d$ by Proposition~(1.5).
So the assertion follows.

\ms\nin
{\bf 4.8.~Proof of Theorem~4 by induction.}
Let $\Phi(A,B,C)$ denote $\Phi_{\S}^0(a)$ in (4.4) where $A=(A_j)$,
$B=(B_k)$, and $A_j,B_k,C$ denote $a_V$ depending on whether
$\dim V=2,1,0$.
Let $a_j,b_k,c$ denote the corresponding divisor classes so that
$\Phi(A,B,C)=\chi(\Y,\E)$ with
$$\E=\O_{\Y}(\msum_jA_ja_j+\msum_kB_kb_k+Cc).$$
We have $\R\pi_*\O_{\Y}=\O_Y$ where $\pi:\Y\to Y$ since $Y$
is nonsingular.
So we first get
$$\h{$\Phi(0,0,C)={C+3\choose 3}$}.$$
Applying the short exact sequence to an exceptional divisor which
is isomorphic to $\widetilde{\P}^2$ and corresponds to some $b_k$,
we get then inductively
$$\h{$\Phi(0,B,C)=\msum_k{B_k\choose 3}+{C+3\choose 3}$},$$
using ${B_k\choose 3}-{B_k-1\choose 3}={B_k-1\choose 2}=
{2-B_k\choose 2}$.
Indeed, the restriction of the exceptional divisor to itself is
the hyperplane section class up to a sign using a general
hyperplane of $\P^3$ passing through the point corresponding
to $b_k$. Then we can use the same argument as above.

We apply the same argument to an exceptional divisor which
is isomorphic to $\P^1\times\P^1$ and corresponds to some $a_j$.
Let $n_j$ be the number of $b_k$ with $k\subset j$.
Here we write $k\subset j$ when there is an inclusion between the
corresponding $V$.
Let ${\bf 1}_j$ be as in the proof of Theorem~(4.4).
We have to calculate
$$\Phi(A,B,C)-\Phi(A-{\bf 1}_j,B,C)=\chi(E,\E|_E).$$
Let $e_1,e_2$ respectively denote the class of $pt\times\P^1$
and $\P^1\times pt$.
Since the restrictions of $a_j$, $b_k\,(k\subset j)$ and $c$ are
respectively
$$(1-n_j)e_1-e_2,\q e_1,\q e_1,$$
the restriction of $\msum_jA_ja_j+\msum_kB_kb_k+Cc$ to
$\P^1\times\P^1$ is
$$\bigl((1-n_j)A_j+\msum_{k\subset j}B_k+C\bigr)e_1-A_je_2.$$
Then
$$\eqalign{\chi(E,\E|_E)&=\bigl((1-n_j)A_j+\msum_{k\subset j}B_k+
C+1\bigr)(1-A_j)\cr
&=\h{$2(n_j-1){A_j\choose 2}-(A_j-1)\bigl(\msum_{k\subset j}B_k+
C+1\bigr)$}.}$$
So we get
$$\h{$\chi(\E)=\msum_j\bigl(2(n_j-1){A_j+1\choose 3}-
{A_j\choose 2}\bigl(\msum_{k\subset j}B_k+C+1\bigr)\bigr)
+\msum_k{B_k\choose 3}+{C+3\choose 3}$}.\leqno(4.8.1)$$
We apply this to $\E$ with
$A_j=2-\lc im_j/d\rc$, $B_j=3-\lc im_j/d\rc$, $C=i-4$,
where $i=d-k$ and $p=3$ in Proposition~(1.5), see (1.4.5).
Then Theorem~4 follows.

\bs\bs
\centerline{\bf 5. Proofs of Theorems~1--4 by HRR}

\bs\nin
{\bf 5.1.~Hirzebruch-Riemann-Roch Theorem.}
For a vector bundle $\E$ of rank $r$ on a smooth complex projective
variety $X$, there are Chern
classes $c_i(\E)\in H^{2i}(X,\Q)$ such that $c_0(\E)=1$,
$c_i(\E)=0$ for $i>r$, and the following facts are well known
(see [17], [19]):

\ms\nin
(a) The total Chern class, the Chern character, and the Todd class
are defined by
$$c(\E)=\msum_i\,c_i(\E),\q ch(\E)=\msum_{1\le i\le r}\exp(x_i),
\q Td(\E)=\mprod_{1\le i\le r}Q(x_i),$$
where $Q(x)=x/(1-\exp(-x))$ and the formal Chern roots $x_i$
satisfy
$$\mprod_{1\le i\le r}\,(1+x_it)=\msum_i\,c_i(\E)t^i.$$
(b) The total Chern class and the Todd class of $X$ are defined by
$$c(X)=c(T_X),\q Td(X)=Td(T_X).$$
(c) By the Hirzebruch-Riemann-Roch theorem [19] we have
$$\chi(\E)=\int_Xch(\E)\,Td(X).\leqno(5.1.1)$$

We will need the following properties of the characteristic classes:

\ms\nin
(d) For a short exact sequence of vector bundles
$0\to\E'\to\E\to\E''\to 0$, we have
$$c(\E)=c(\E')\,c(\E''),\q ch(\E)=ch(\E')+ch(\E''),\q
Td(\E)=Td(\E')\,Td(\E'').\leqno(5.1.2)$$
(e) For the tensor product of two vector bundles $\E,\F$ we have
$$ch(\E\h{$\otimes$}\F)=ch(\E)\,ch(\F).\leqno(5.1.3)$$
(f) For the exterior product we have
$$\h{$\sum_ic_i\bigl(\mwedge^p\E\bigr)t^i=\prod_{i_1<\cdots<i_p}
(1+(x_{i_1}+\cdots x_{i_p})t)$}.\leqno(5.1.4)$$
(g) For the dual vector bundle $\E^{\vee}$, we have
$$c_i(\E^{\vee})=(-1)^ic_i(\E).\leqno(5.1.5)$$

\ms\nin
{\bf 5.2.~Remarks.} (i) The function $Q(x)$ has the expansion
$$Q(x)=1+{1\over 2}x+\sum_{k=1}^{\infty}(-1)^{k-1}{B_k\over(2k)!}
x^{2k},\leqno(5.2.1)$$
where the $B_k$ are the Bernoulli numbers, see [17], Ex.~3.2.4.
The first few terms of $B_k$ are
$${1\over 6},\,{1\over 30},\,{1\over 42},\,{1\over 30},\,{5\over 66},
\,\dots$$

\ms
(ii) Using $c_i=c_i(\E)$ and $r={\rm rank}\,\E$, we have the
expansions (see [17], Ex.~3.2.3--4)
$$\eqalign{ch(\E)&=r+c_1+{1\over 2}(c_1^2-2c_2)+{1\over 6}
(c_1^3-3c_1c_2+c_3)+\cdots,\cr
Td(\E)&=1+{1\over 2}c_1+{1\over 12}(c_1^2+c_2)+{1\over 24}
(c_1c_2)+\cdots.}\leqno(5.2.2)$$

\ms
(iii) By (5.1.2), $c(\E),ch(\E),Td(\E)$ are extended to
well-defined morphisms
$$c(\E),ch(\E),Td(\E):K^0(X)\to H^{\ssb}(X,\Q),\leqno(5.2.3)$$
where the source is the Grothendieck group of vector bundles on $X$.
Note that the initial term of $ch(\E)$ is the virtual rank of $\E$,
and the latter does not appear in $c(\E),Td(\E)$.

\ms
(iv) For $n=\dim X$ we have
$$\int_X c_n(X)=\chi(X,\C),\q \int_X Td(X)_n=\chi(X,\O_X),
\leqno(5.2.4)$$
where $\chi(X,\C)$ is the topological Euler characteristic of $X$.
For the first assertion, see e.g.\ [17], Ex.~8.1.12.
The second assertion follows from the Hirzebruch-Riemann-Roch
theorem (5.1.1) applied to $\E=\O_X$ where $c_i(\E)=0$ for $i>0$.

\ms\nin
{\bf 5.3.~Combinatorial description of the cohomology.}
Let $D$ be an essential central hyperplane arrangement.
In the notation of (1.1) we apply the construction in (2.1) to
$$\S:=\S(D)^{\nnc},$$
and not to $\S(D)$ as in [4], [5].
(This simplifies some arguments in loc.~cit.\ considerably.)
Note that $\g(V)\ge 2$ for $V\in\S$.
By C.~De Concini and C.~Procesi [8] the cohomology ring of $\Y$
in (2.1) is described by using only the combinatorial data as
follows:

Let $\Q[e_V]_{V\in\S}$ be the polynomial ring with independent
variables $e_V$ for $V\in\S$.
There is an isomorphism
$$\Q[e_V]_{V\in\S}/I_{\S}\simto H^{\ssb}(\Y,\Q),\leqno(5.3.1)$$
sending $e_V$ to $[E_V]$ for $V\ne 0$ and $e_0$ to $-[E_0]$,
where $E_0$ is the total (or proper) transform of a general
hyperplane which was denoted by $\H$.
Moreover, the ideal $I_{\S}$ is generated by
$$R_{V,W}=\cases{e_Ve_W&if $V,W$ are incomparable,\cr
e_V\e_W^{\,\g(W)-\g(V)}\raise12pt\h{$\,$}&if $W\ssn V$,\cr
\e_W^{\,\g(W)}&if $V=\C^n$,}\leqno(5.3.2)$$
where $\e_W:=\sum_{W'\subset W}e_{W'}$ and $\g(V):=\codim\,V$.
Here $V,W,W'\in\S$ except for the third case where $V=\C^n$.
Note that $\S$ is stable by intersection so that a {\it nested}
subset of $\S$ in the sense of [8] is always linearly ordered
by the inclusion relation.

For $V\in\S(D)\setminus\S(D)^{\nnc}$, let $\P(V)^{\sim}$ denote
the proper transform of $\P(V)$ in $\Y$.
(Here the notation $\P(V)^{\S^V}$ in Section 2 cannot be used
since $V\notin\S:=\S(D)^{\nnc}$.)
Then the cohomology class $e_V$ of $\P(V)^{\sim}$ is given by
$\prod_{D_j\supset V}e_{D_j}$ since $\P(V)^{\sim}$ is the
intersection of $\P(D_j)^{\sim}$ with $D_j\supset V$.
For $V=D_j$, we have by calculating the total transform of $D_j$
$$e_{D_j}+\msum_{W\in\S,W\subset{D_j}}\,e_W=0,\leqno(5.3.3)$$
since $e_0$ corresponds to $-\H$.
(This is similar to (5.3.2) for $\S=\S(D)$ although $\Y$ is
different.)

\ms\nin
{\bf 5.4.~Calculation of Chern classes.}
In our case the Chern classes of $\Y$ are expressed by applying
inductively the formula for the Chern classes under the blow-up
([17], 15.4).
By [4] we have $c(\Y)=\mprod_{V\in\S}F_V$, under the isomorphism
(5.3.1), with
$$F_V=\cases{(1-\ee_V)^{-\g(V)}(1+e_V)(1-\e_V)^{\g(V)}&
if $V\ne 0$,\cr (1-e_0)^n&if $V=0$,}\leqno(5.4.1)$$
where $\e_V:=\sum_{W\subset V}e_W$, $\ee_V:=\e_V-e_V$, and
$\g(V):=\codim\,V$.
Using the Grothendieck group as in (5.2.3), (5.4.1)
implies that $Td(\Y)=\mprod_{V\in\S}G_V$ with
$$G_V=\cases{Q(-\ee_V)^{-\g(V)}Q(e_V)Q(-\e_V)^{\g(V)}&
if $V\ne 0$,\cr Q(-e_0)^n&if $V=0$,}\leqno(5.4.2)$$
So the Chern classes and the Todd class of $\Y$ are
expressed by using only the combinatorial data via (5.3.1).

Set $\S'=\S(D)^{\nnc}\cup\{D_i\}$ where the $D_i$ are the irreducible
components of $D$.
The proper transform $\P(D_i)^{\sim}$ of $\P(D_i)$ in $\Y$ will be
denoted by $E_{D_i}$ .
By (5.3.3) its cohomology class $e_{D_i}$ is given by
$$e_{D_i}=-\msum_{W\in\S,\,W\subset D_i}\,e_W\in\Q[e_W]_{W\in\S}/
I_{\S}.$$
We have a short exact sequence
$$0\to\Omega_{\Y}^1\to\Omega_{\Y}^1(\log\Z)\to\mopl_{V\in
\S'\setminus\{0\}}\,\O_{E_V}\to0.$$
Using (5.2.3), we get then
$$c(\Omega_{\Y}^1(\log\Z))=c(\Omega_{\Y}^1)\mprod_V\,c(\O_{E_V})
=c(\Omega_{\Y}^1)\mprod_V\,c(\O_{\Y}(-E_V))^{-1}.\leqno(5.4.3)$$
Moreover, the Chern classes of $\Omega_{\Y}^p(\log\Z)
=\mwedge^p\Omega_{\Y}^1(\log\Z)$ for $p>1$ are expressed by using
those of $\Omega_{\Y}^1(\log\Z)$ by (5.1.4).
(However, it is not easy to write down the universal polynomials
explicitly.)

\ms\nin
{\bf 5.5.~Proofs of Theorems~1--2 by HRR.}
We calculate the right-hand side of (1.5.1) in Proposition~(1.5)
by applying the Hirzebruch-Riemann-Roch theorem~(5.1.1) to
$$\E_k=\Omega_{\Y}^{p}(\log\Z)\otimes_{\O}\O_{\Y}
(-k\H+\msum_j\,\lf km_j/d\rf E_j),\leqno(5.5.1)$$
where the $m_j$ are given by $\mu(V)$ if $E_j$ in Proposition~(1.5)
is $E_V$ in (2.1).
For $\O_{\Y}(D_k)$ with
$$D_k=-k\H+\msum_j\,\lf km_j/d\rf E_j,$$
we have $c(\O_{\Y}(D_k))=1+[D_k]$.
Then we can apply (5.1.3) to calculate $ch(\E_k)$,
and $\chi(\E_k)$ depends only on the combinatorial data using
the assertions in (5.4) together with the Hirzebruch-Riemann-Roch
theorem (5.1.1).
In the reduced case the $n_{f,\a}$ depend only on the
combinatorial data as in Theorem~2.
Moreover, if $p=0$ or $p=n-1$, then $n_{f,\a}$ for
$\a\in(0,1]\cup(n-1,n)$ depends only on the weak equivalence class
using (1.4.5) for $p=n-1$ since $\Omega_{\Y}^0(\log\Z)=\O_{\Y}$
for $p=0$.
So the assertion follows.

\ms\nin
{\bf 5.6.~Remark.}
We can prove Theorem~(4.4) by using (5.3--4), and this is enough
for the proofs of Theorems~1--2 as is shown in (4.5).

In the following, we illustrate how to calculate $n_{f,\a}$ using
this method.

\ms\nin
{\bf 5.7.~Proof of Theorem~3 by HRR.}
Let $a_i$ denote the $e_V$ mod $I_{\S}$ in (5.3) for $V\in\S^{(2)}$
(see (1.1.3)), i.e. the $a_i$ correspond to the points of
$\P(D^{\nnc})$.
Set $c=e_0$.
We have the relations
$$a_ia_j=0\,(i\ne j),\,\,a_ic=0,\,\,a_i^2=-c^2,\,\,c^3=0,$$
using $(a_i+c)^2=0$, etc. in (5.3.2).
Let $F'_i$ denote $F_V$ for $V$ corresponding to $a_i$. Then
$$F'_i=(1-c)^{-2}(1+a_i)(1-c-a_i)^2=1-a_i+c^2,$$
and $F_0=(1-c)^3$.
Set $\nu^{(2)}=\sum_{m\ge 3}\nu^{(2)}_m$ with $\nu^{(2)}_m$ as in
(1.1.3).
Since $c(\Y)=F_0\mprod_iF'_i$, we get
$$c(\Y)=1-(\msum a_i+3c)+(\nu^{(2)}+3)c^2,\q
Td(\Y)=1-{1\over 2}(\msum a_i+3c)+c^2,$$
using (5.2.2).
Note that $\Omega_{\Y}^2=\O_{\Y}(-3\H+\msum_i E_i)$ where the
$E_i$ correspond to $a_i$.
So the Hirzebruch-Riemann-Roch theorem for a line bundle
coincides with the Riemann-Roch theorem for surfaces, and
the argument is the same as in (3.1) if $p=2$ or $0$.

In (3.1), the assertion for $p=1$ is reduced to the other cases
using the relation with $\chi(U)$.
However, it should be possible to prove it by using the
Hirzebruch-Riemann-Roch theorem for vector bundles.
We apply this to
$$\E_k=\Omega_{\Y}^1(\log\Z)\otimes\O_{\Y}(D_k),$$
with $D_k=-k\H+\msum_j\,\lf km_j/d\rf E_j$.
By the calculation of $c(\Y)$ together with (5.1.5) we have
$$c(\Omega_{\Y}^1)=1+\msum a_i+3c+(\nu^{(2)}+3)c^2.$$
So we get by (5.4.3)
$$c(\Omega_{\Y}^1(\log\Z))=(1+\msum a_i+3c+(\nu^{(2)}+3)c^2)\,
\mprod_i\,(1-a_i)^{-1}\,\mprod_j\,(1-a'_j)^{-1}.$$
Here $a'_j:=-\sum_{i\subset j}a_i-c$ which corresponds to the
proper transform of an irreducible component $\P(D_j)$ of $\P(D)$,
and we write $i\subset j$ if the corresponding subspaces of $\P^2$
have such an inclusion relation.
(Note that $-c$ corresponds to $\H$.)
We have
$$\h{$\mprod_j\,(1+\msum_{i\subset j}a_i+c)=1+\msum_i\,m_ia_i+dc
+\bigl({d\choose 2}-\msum_i\,{m_i\choose 2}\bigr)c^2$,}$$
and $c(\Omega_{\Y}^1(\log\Z))$ is equal to
$$1+\msum_i(2-m_i)a_i+(3-d)c-\h{$1\over 2$}
\bigl((d^2-5d+2\nu^{(2)}+6)-\msum_i(m_i^2-3m_i+4)\bigr)c^2.$$
Then $Td(\Y)$, $\h{$1\over 2$}ch(\Omega_{\Y}^1(\log\Z))$,
$ch(\O_{\Y}(D_k))$ are respectively
$$\eqalign{&1-\h{$1\over 2$}(\msum_ia_i+3c)+c^2,\cr
&1-\h{$1\over 2$}\bigl(\msum_i(m_i-2)a_i+(d-3)c\bigr)
+\h{$1\over 4$}\bigl(\msum_im_i-d-2\nu^{(2)}+3\bigr)c^2,\cr
&1+\bigl(\msum_i\lf km_i/d\rf a_i+kc\bigr)
-\h{$1\over 2$}\bigl(\msum_i\lf km_i/d\rf^2-k^2\bigr)c^2.}$$
Calculating the degree 2 part of the multiplication of these three,
we get the right-hand side of the second equation divided
by $-2$ in Theorem~3 where $i=d-k$.

\ms\nin
{\bf 5.8.~Proof of Theorem~4 by HRR.}
Let $a_j\,(j\in I')$, $b_k\,(k\in I'')$, $c$ denote the $e_V$ mod
$I_{\S}$ in (5.3) for $V\in\S^{(i)}$ (see (1.1.3)) depending on
whether $i=2$ or 3 or 4.
We will write $k\subset j$ if there is an inclusion relation
between the corresponding $V$.
Let $n_j\,(j\in I')$ be the number of $b_k\,(k\in I'')$ with
$k\subset j$.
We have the relations
$$\eqalign{&a_ia_j=b_kb_l=a_jb_k^2=a_jc^2=b_kc=0\,(i\ne j,\,k\ne l),
\,a_jb_k=0\,(k\not\subset j),\cr &a_jb_k=-a_jc\,(k\subset j),\,a_j^3
=-2(n_j-1)c^3,\,a_j^2c=b_k^3=-c^3,\,a_j^2b_k=c^3\,(k\subset j),}$$
using $(a_j+\msum_{k\subset j}b_k+c)^2=0$, $(b_k+c)^3=0$,
$a_j(b_k+c)=0\,(k\subset j)$, $a_jc^2=0$, see (5.3.2).
By the same argument as in (4.8) it is enough to calculate
$\chi(\E)$ for a line bundle $\E$ with
$$c_1(\E)=\msum_jA_ja_j+\msum_kB_kb_k+Cc\,\,\,(A_j,B_k,C\in\ZZ).
\leqno(5.8.1)$$

Let $F'_j,\,F''_k$ denote $F_V$ if $V$ corresponds to $a_j,b_k$
respectively. Then
$$F'_j=1-a_j-a_j^2+2(n_j-1)a_jc-2(n_j-1)c^3,\q
F''_k=1-2b_k-2c^3.$$
Since $c(\Y)=F_0\mprod_jF'_j\mprod_kF''_k$ with $F_0=(1-c)^4$,
we get
$$c_1(\Y)=-\msum_ja_j-\msum_k2b_k-4c,\q
c_2(\Y)=\msum_j(2a_jc-a_j^2)+6c^2,$$
where $c_3(\Y)$ is the topological Euler characteristic $\chi(\Y)$
multiplied by $-c^3$, see (5.2.4).
This gives $Td(\Y)$ using the expansion of the Todd class in (5.2.2)
(where $c_3$ does not appear so that $c_3(\Y)$ is not needed).
Thus we get
$$\eqalign{Td(\Y)_1&=-\msum_j{a_j\over 2}-\msum_kb_k-2c,\cr
Td(\Y)_2&=\msum_j{(5-2n_j)a_jc\over 6}+\msum_k{b_k^2\over 3}+
{11c^2\over 6},}$$
where $Td(\Y)_3=-c^3$, see (5.2.4).
We apply this theorem to the line bundle $\E$ in (5.8.1).
Then we get
$$\h{$\chi(\E)=\msum_j\bigl(2(n_j-1){A_j+1\choose 3}-
{A_j\choose 2}\bigl(\msum_{k\subset j}B_k-C+1\bigr)\bigr)
+\msum_k{B_k\choose 3}-{C-1\choose 3}$},\leqno(5.8.2)$$
which is compatible with (4.8.1) where $C$ corresponds to $-C$.
We apply this to $\E$ with
$A_j=2-\lc im_j/d\rc$, $B_j=3-\lc im_j/d\rc$, $C=4-i$,
where $i=d-k$ and $p=3$ in Proposition~(1.5), see (1.4.5).
Then Theorem~4 follows.

\vfill\eject
\centerline{{\bf References}}

\ms
{\mfont
\item{[1]}
Beilinson, A., Bernstein, J.\ and Deligne, P., Faisceaux pervers,
Ast\'erisque, vol. 100, Soc. Math. France, Paris, 1982.

\item{[2]}
Brieskorn, E., Sur les groupes de tresses
[d'apr\`es V.I. Arnold], S\'eminaire Bourbaki, 24\`eme ann\'ee
(1971/1972), Exp. No. 401, Lect. Notes in Math.
Vol.~317, Springer, Berlin, 1973, pp. 21--44.

\item{[3]}
Budur, N., On Hodge spectrum and multiplier ideals,
Math. Ann.\ 327 (2003), 257--270.

\item{[4]}
Budur, N., Jumping numbers of hyperplane arrangements,
arXiv:0802.0878 to appear in Comm.\ Algebra.

\item{[5]}
Budur, N., Hodge spectrum of hyperplane arrangements,
arXiv:0809.3443 (unpublished).

\item{[6]}
Budur, N.\ and Saito, M., Multiplier ideals, $V $-filtration, and
spectrum, J.\ Alg.\ Geom.\ 14 (2005), 269--282.

\item{[7]}
Cohen, D.C.\ and Suciu, A., On Milnor fibrations of arrangements,
J.\ London Math.\ Soc.\ 51 (1995), 105--119.

\item{[8]}
De Concini, C.\ and Procesi, C., Wonderful models of subspace
arrangements, Selecta Math.\ (N.S.) 1 (1995), 459--494.

\item{[9]}
Deligne, P.,
Equations Diff\'erentielles\`a Points Singuliers R\'eguliers,
Lect.\ Notes in Math.\ vol.\ 163, Springer, Berlin, 1970.

\item{[10]}
Deligne, P., Th\'eorie de Hodge II, Publ.\ Math.\ IHES, 40 (1971),
5--58.

\item{[11]}
Deligne, P., Le formalisme des cycles \'evanescents, in SGA7 XIII
and XIV, Lect.\ Notes in Math.\ 340, Springer, Berlin, 1973,
pp.\ 82--115 and 116--164.

\item{[12]}
Dimca, A., Singularities and Topology of Hypersurfaces,
Universitext, Springer, Berlin, 1992.

\item{[13]}
Dimca, A.\ and Saito, M., A generalization of Griffiths's theorem
on rational integrals, Duke Math.\ J.\ 135 (2006), 303--326.

\item{[14]}
Ein, L., Lazarsfeld, R., Smith, K.E.\ and Varolin, D., Jumping
coefficients of multiplier ideals, Duke Math.\ J.\ 123 (2004),
469--506.

\item{[15]}
Esnault, H., Fibre de Milnor d'un c\^one sur une courbe plane
singuli\`ere, Inv.\ Math.\ 68 (1982), 477--496.

\item{[16]}
Esnault, H.\ and Viehweg, E., Rev\^etements cycliques,
in Algebraic threefolds (Varenna, 1981), Lecture Notes in Math.,
947, Springer, Berlin-New York, 1982, pp.~241--250.

\item{[17]}
Fulton, W., Intersection Theory, Springer, Berlin, 1984.

\item{[18]}
Hartshorne, R., Algebraic Geometry, Springer, Berlin, 1977.

\item{[19]}
Hirzebruch, F., Topological method in algebraic geometry,
Springer, Berlin, 1966.

\item{[20]}
Jouanolou, J.-P., Cohomologie de quelques sch\'emas classiques
et th\'eorie cohomologique des classes de Chern, in SGA5, Lecture
Notes in Math.\ 589, Springer, Berlin, 1977, pp.~282--350.

\item{[21]}
Lazarsfeld, R.,
Positivity in algebraic geometry II,
Ergebnisse der Mathematik und ihrer Grenzgebiete.\ 3.~Folge,
A series of Modern Surveys in Mathematics, Vol.~49,
Springer-Verlag, Berlin, 2004.

\item{[22]}
{\Mustata}, M., Multiplier ideals of hyperplane arrangements,
Trans.\ Amer.\ Math.\ Soc.\ 358 (2006), 5015--5023.

\item{[23]}
Nadel, A.M., Multiplier ideal sheaves and K\"ahler-Einstein
metrics of positive scalar curvature, Ann.\ Math.\ 132 (1990),
549--596.

\item{[24]}
Orlik, P.\ and Solomon, L., Combinatorics and topology of
complements of hyperplanes, Inv.\ Math.\ 56 (1980), 167--189.

\item{[25]}
Rybnikov, G., On the fundamental group of the complement of a
complex hyperplane arrangement (math.AG/9805056).

\item{[26]}
Saito, M., Mixed Hodge modules, Publ.\ RIMS, Kyoto Univ.\ 26
(1990), 221--333.

\item{[27]}
Saito, M., Multiplier ideals, $b$-function, and spectrum of a
hypersurface singularity, Compos.\ Math.\ 143 (2007) 1050--1068.

\item{[28]}
Saito, M., Bernstein-Sato polynomials of hyperplane arrangements
(math.AG/0602527).

\item{[29]}
Saito, M., Jumping coefficients and spectrum of a hyperplane
arrangement (unpublished manuscript, 2007).

\item{[30]}
Schechtman, V., Terao, H.\ and Varchenko, A.,
Local systems over complements of hyperplanes and the Kac-Kazhdan
conditions for singular vectors, J.\ Pure Appl.\ Algebra 100 (1995),
93--102.

\item{[31]}
Steenbrink, J.H.M., The spectrum of hypersurface singularity,
Ast\'erisque 179--180 (1989), 163--184.

\item{[32]}
Teitler, Z., A note on {\Mustata}'s computation of multiplier
ideals of hyperplane arrangements, Proc.\ Amer.\ Math.\ Soc.\ 136
(2008), 1575--1579.

\item{[33]}
Walther, U., Bernstein-Sato polynomial versus cohomology of the
Milnor fiber for generic hyperplane arrangements,
Compos.\ Math.\ 141 (2005), 121--145.

{\sfont\baselineskip=10pt
\ms
Department of Mathematics, The University of Notre Dame, IN 46556,
USA

e-mail: nbudur@nd.edu

\ss
RIMS Kyoto University, Kyoto 606-8502 Japan

e-mail: msaito@kurims.kyoto-u.ac.jp

\ms\vers}}
\bye